\documentclass[11pt,a4paper,leqno]{amsart}
\usepackage{amsmath,amsfonts}
\usepackage{amssymb}
\usepackage{upref}
\usepackage[mathcal]{eucal}
\usepackage{graphicx}

\date{\today}

\begin{document}
\title{On sequence  groups}
\author{Zbigniew Lipi{\'n}ski}
\address{University of Opole\\
 Opole, Poland}
\email{zlipinski@uni.opole.pl}

\author{Maciej P. Wojtkowski}
\address{ The University of Arizona \\
 Tucson, Arizona 85721, USA}
\email{maciejw@math.arizona.edu}

\date{\today}

            \theoremstyle{plain}
\newtheorem{lemma}{Lemma}%[section]
\newtheorem{proposition}[lemma]{Proposition}
\newtheorem{theorem}[lemma]{Theorem}
\newtheorem{corollary}[lemma]{Corollary}
\newtheorem{fact}[lemma]{Fact}
   \theoremstyle{definition}
\newtheorem{definition}{Definition}
    \theoremstyle{example}
\newtheorem{example}{Example}
    \theoremstyle{remark}
\newtheorem{remark}{Remark}

\newcommand{\macierz}[4]{\scriptsize{\left[\begin{array}{cc} #1 & #2 \\
        #3 & #4 \end{array}\right]}}

\begin{abstract}
Linear second order recursive sequences with arbitrary initial conditions
  are studied. For sequences with the same parameters a ring and a group
  is attached, and isomorphisms and homomorphisms are established
  for related parameters.
  In the group, called the {\it sequence group},
  sequences are identified if they differ by a scalar factor,
  but not if they differ by a shift, which  is the case for the
  Laxton group.

  Prime divisors of sequences
  are studied with the help of the sequence group $\mod p$, which is always
  cyclic of order $p\pm 1$.

  Even and odd numbered subsequences are given independent status through
  the introduction of one rational parameter in place of two integer
  parameters. This step brings significant simplifications in the
  algebra.

  All elements of finite order in Laxton groups and sequence
  groups are described effectively.

 A necessary condition is established for  a prime $p$
 to be a divisor of a sequence: {\it the norm (determinant) of
   the respective
    element of the ring must be a square $\mod p$}.
    This leads to an uppers estimate of the set of  divisors
    by a set of  prime density $1/2$. Numerical experiments show that
    the actual density is typically close to $0.35$.

  A conjecture is formulated that the sets of prime divisors of the
  even and odd numbered elements are independent for a large family
  of parameters.
\end{abstract}

\maketitle

\section{Introduction}

We present some new results about linear recursive sequences
of order $2$. There is a vast literature of the subject.
The book by
Everest,  van der Poorten, Shparlinski and  Ward, \cite{E-P-Sh-W},
gives a broad panorama of problems and results with an
exhaustive bibliography.  Another source is
the book of Williams, \cite{Wi}, where the  themes developed
in our work are brought into focus.
We will refer only to papers which are most relevant to our
work. Our study is rooted in the work of Laxton, \cite{L1, L2}
who associated an abelian group structure to recursive sequences,
to study the sets of their prime divisors.

For  rational  $ Q\neq 0$ and  $T \neq 0$  let $D= D_{T,Q} =
\left[ \begin{array}{cc} 0 &     -Q \\
    1  &   T \end{array} \right]$
be the matrix defining recursive sequences $\{x_n\}_{n\in\mathbb Z}$ by the formula
\[
[x_n \  x_{n+1}]=[x_0 \  x_{1}]D^n.
\]
We choose to consider all rational initial conditions
$x_0,x_1$. Non-integer values are actually unavoidable even if
we consider integer parameter pairs $(T,Q)$ and integer initial values.
Indeed powers of $Q$ appear in the
denominators of the elements with negative indices.
In general the denominators
are similarly ``poor'',
while the numerators are ``rich'', in prime divisors.
More precisely for every sequence
the denominators have only finitely many prime divisors.
This allows the use of the modular
arithmetic with rational numbers. A fully reduced simple fraction
$\frac{a}{b} = 0 \mod p$, for a prime $p$ not dividing $b$,
if $a$ is divisible by $p$.
It works well if we avoid
the prime divisors of the denominators, which are only finitely
many for a specific sequence.

 We propose to  separate
  the even and odd terms of a recursive sequence as recursive sequences
  in their own right. More precisely
one  sequence $\{x_n\}_{n\in\mathbb Z}$ for the parameter pair
$(T,Q)$ is split into two  sequences
$\{y_k =Q^{1-k}x_{2k}\}_{k\in\mathbb Z}$ and $\{z_k =Q^{1-k}x_{2k-1}\}_{k\in\mathbb Z}$
for the parameters $(t,1), t = T^2/Q-2$,
essentially the even and odd numbered elements
of the original sequence.
The idea to allow the rational parameter $t$,
  and achieve $Q=\det D = 1$, appeared first in \cite{Wo},
  and it lead to substantial simplifications. In this paper
  we find further applications of this method. We mostly use
  the one parameter $t$ language, due to its simplicity,
  and translate the results into the language of two parameters
  $(T,Q)$.

  The interest in Laxton groups was recently revived
  by Aoki and Kida, \cite{A-K}, and Suwa, \cite{Suwa1, Suwa2}.

  In Section 2 we introduce a ring
  of $2\times 2$ matrices associated with recursive
  sequences. In general this ring is  a quadratic
  field extension of $\mathbb Q$. Similar construction appeared
  in the work of Hall, \cite{H1}. We use it to compare sequences for
  the parameters $t$ and $t_r = C_r(t)$, where
  $C_r(t) = tr \ D^r_t$ are the Chebyshev polynomials of the first kind.
  In Proposition \ref{Proposition4}  we show that the sequences for $t$
  can be disintegrated into $r$ sequences for $t_r$.
  It brings forward the concept of a {\it primitive}
  parameter $t$ which is not equal to $C_r(u)$
  for any prime $r$ and  rational $u$.

  In Section 3 we discuss the phenomenon of {\it recombination}:
  two sequences for a parameter pair $(T,Q)$
  become sequences for another
  parameter pair
  $(\widehat T, \widehat Q)$, after essentially the exchange of odd terms.
The two parameter pairs are rigidly connected, and we call them {\it twins}.
For example the twin of the Fibonacci parameters $(1,-1)$ is $(5,5)$.
The phenomenon of twins was discovered in \cite{Wo} in a more
restricted setting of special sequences.

We proceed looking for rational $t\neq  a$ and a prime $r$
such that
$C_r(t) = C_r(a)$, which is
essential for extending the recombination to other sets
of parameters. It turns out to be   exceedingly rare,
 limited to  $r=2$ (the {\it twin} case), $r=4$ (the {\it circular} case) and
 $r=3$ (the {\it cubic} case). The parameter $a$ will be called
an  {\it associate} parameter of $t$.

% The first case  is that of the twins
%$t\leftrightarrow -t$ for $r=2$.

In Section 4 we introduce the {\it sequence group} $\mathcal L(t)$,
which is an  extraction of the multiplicative
structure of the ring up to scalar factors.
Further we consider $\mathcal L_p(t)$,  the {\it sequence group $\mod p$},
which turns out to be always cyclic of order $p\pm 1$ (Theorem \ref{Theorem16}).

In Section 5 we arrive at the Laxton group $\mathcal G(t)$, \cite{L1},
which is the  quotient of the sequence group by the
cyclic subgroup generated by the element $D$.
While the finite order elements in $\mathcal G$
are of significant interest, the sequence groups $\mathcal L$ provide
a good environment for calculations.

We are able to list all torsion elements in sequence groups
and Laxton groups for all
parameters. It was attempted in \cite{L2}, we give a more
explicit description.
For primitive parameters $t$ there are only $3$ finite order
elements in the Laxton group under additional assumption of
primitivity, namely that not only $t$ is primitive but also
its associate parameters are primitive. This stronger condition
is called {\it twin}, {\it circular} or  {\it cubic} primitivity.

%except in the circular and cubic cases. The circular case
%is especially involved, it requires the introduction of
%the refinement of primitivity to {\it circular primitivity}.

The special sequences which are of finite order in the circular
and cubic cases were studied in detail in a remarkable paper  of
Ballot, \cite{B}, in the language of the parameter pairs $(T,Q)$.

In Section 6 we briefly describe the structure of all the torsion subgroups
of the Laxton groups. However we must admit that we did not find
any application for this information.

In Section 7 we proceed with the study of the sets of prime divisors
of recursive sequences, mostly in the language of one rational
parameter $t$.
We reformulate the old theorem of Hall, \cite{H2}, into
the language of $\mathcal L_p(t)$ (Proposition \ref{Proposition29}).

In his paper, \cite{L1},  Laxton states that
{\it ``One problem of perennial interest is that of prime divisors of a recurrence. ...  The divisor problem is the chief interest in this article.''}

The study of the sets of prime divisors for recursive sequences
for arbitrary initial conditions was the subject of only few papers,
Laxton was preceded by Ward, \cite{Wa},
where he proves that for most sequences this set
is infinite.  The exclusions are only in parameters, and not the
initial conditions. In our one parameter language they correspond to
$t\neq 0,\pm 1, \pm 2$. This result can be also traced back to the paper
of Polya, \cite{Polya}, \cite{Leeuwen}.
The results of Stephen, \cite{St},
and Moree and Stevenhagen, \cite{M-S},
suggest that the sets have positive prime density.
The state of the art is presented comprehensively by Moree
in  his survey paper \cite{M}, Chapter 8.4.

In the final Section 8 we formulate a conjecture
that  for a large family of sequences,
for the parameter pair $(T,Q)$, the sets of divisors
of even and odd terms are {\it independent}, i.e.,
the prime density of their intersection is the product
of their prime densities.  This conjecture has two
motivations.
One is the simple criterion (Theorem \ref{Theorem31}) that
for most sequences if an odd prime $p$ is a divisor
then the determinant (or the norm) of the sequence is
a square $\mod p$. This implies that for most
sequences of one  parameter $t$ the set of divisors
is contained in a set of prime density $1/ 2$.
As a consequence we arrive at two independent sets
for the even and odd terms of the sequence for
the respective parameter pair $(T,Q)$.

The other motivation comes from numerical experiments
which give a tentative confirmation of the conjecture.

\section{The ring of recursive sequences}
%$\in \mathbb Z[T,Q^{-1}]$
%then all the elements of the sequence
%belong  to $\mathbb Z[T,Q^{-1}]$. The reason for restricting
%the coefficients to this ring, rather than the field of all rationals,
%is that it lends itself readily to the homomorphic factorization
%onto $\mathbb F_p$, with the exclusion of only finitely many primes $p$
%that divide the denominators of $T$ or $Q^{-1}$.

Let $\mathcal R = \mathcal R(T,Q)$ be the ring of
$2\times 2$ matrices with rational
entries which commute with $D= D_{T,Q} =
\left[ \begin{array}{cc} 0 &     -Q \\
    1  &   T \end{array} \right]$. In particular $D \in \mathcal R$.
\begin{lemma}
 \label{Lemma1}
  \[
  \mathcal R = \left\{X\ | \
X=
\left[ \begin{array}{cc} -Qx_{-1} &     -Qx_{0} \\
    x_{0} &   x_{1} \end{array} \right], \ x_1 = Tx_0-Qx_{-1},
\ x_0,x_1 \in \mathbb Q\right\}
\]
\end{lemma}
For ease of notation we will denote an element $X \in \mathcal R$
by the second row $[x_0 \   x_1]$.
We put the recursive
sequences $\{x_n\}_{n\in\mathbb Z}$ into
$1-1$ correspondence with the elements of $\mathcal R$ by the formula
$X= [x_0\ x_1]$. Alternatively we associate with $\{x_n\}_{n\in\mathbb Z}$
the sequence of matrices $\{XD^n\}_{n\in\mathbb Z}$ in $\mathcal R$, where
\[
XD^n=
\left[ \begin{array}{cc} -Qx_{n-1} &     -Qx_{n} \\
    x_{n} &   x_{n+1} \end{array} \right].
\]

%, with initial data $x_0,x_1$. In other words we can naturally
%identify $\mathcal R$ with $\mathbb Z[T,Q^{-1}]\times \mathbb Z[T,Q^{-1}]$.
 The determinant of the matrices in
$\mathcal R$ will play an important role.
    For $[x_0  \  x_1]=X\in\mathcal R$
    \[
    \det \ X = x_1^2-Tx_1x_0+Qx_0^2,
      \]
      As long as the discriminant $\Delta = T^2-4Q$ is not a rational square
      there are no zero divisors in the ring, and it is actually
      isomorphic to the quadratic field $\mathbb Q\left(\sqrt{\Delta}\right)$.
      Under this isomorphism the determinant becomes  the norm
      of a field element.
The canonical field automorphism of $\mathbb Q\left(\sqrt{\Delta}\right)$
 translates into $\mathcal R(T,Q)$ as
$X \mapsto (\det X) X^{-1}$. If $X$ represents a recursive sequence
$\{x_n\}_{n\in\mathbb Z}$ then $(\det X) X^{-1}$ represents
the sequence $\{-Q^nx_{-n}\}_{n\in\mathbb Z}$.

By direct calculation we obtain that the traces $\widehat x_n = tr\ XD^n$
form a sequence in  $\mathcal R(T,Q)$
equal to $\widehat X = CX$ where the sequence $C = [2 \ T]$.
In a somewhat different formulation, Laxton, \cite{L1}, called the sequences
$X$ and $\widehat X$ {\it polar} sequences. Since $C^2 = -\Delta I$,
where $\Delta = T^2-4Q$, we have that polarity is a symmetric relation
defined by the property that  polar sequences have squares
differing only by a scalar, namely $\widehat X^2 = -\Delta X^2$.
The polar sequences are connected
in the following equation.
\[
\widehat x_n^2 - 4(\det X)Q^n = \Delta\ x_n^2
\]
Laxton, \cite{L1} (Theorem 4.3), discovered for any sequence
$X$ the following factorization of elements of $X^2$.
\begin{proposition}
 \label{Proposition6}
  For any  $[x_0 \ x_1] = X \in \mathcal R(T,Q)$ we have
  \[
   X^2 =\left[ \begin{array}{cc} -Qx_{-1} &     -Qx_0 \\
      x_0 &   x_1 \end{array} \right]^2=
  \left[ \begin{array}{cc}
      x_0(x_1-Qx_{-1}) &   x_1^2-Qx_0^2 \end{array} \right],
  \]
\end{proposition}
Hence the even numbered elements in $X^2$
factor into elements of $X$ and $\widehat X = CX  = [x_{1}-Qx_{-1} \ x_2-Qx_0]$.
Odd numbered elements also factor into sequences from $\mathcal R(T,Q)$
in the case when $Q$ is a square.

  We introduce now the crucial step in our paper, namely the separation
  of even and odd terms of a recursive sequence as recursive sequences
  in their own right.
%  \begin{proposition}
%    For any pair of rational nonzero parameters $(T,Q)$  and $t = T^2Q^{-1}-2$
%    there is a  canonical ring isomorphism $\Phi: \mathcal R(T,Q)
%    \to \mathcal R(t,1)$,
%    $\mathcal L(T,Q)$ and $\mathcal L(t,1)$,
%    defined  by the conjugation $\Phi(X)=  A^{-1}XA  \in \mathcal R(t,1)$,
%    where
%$
%A =\left[ \begin{array}{cc} Q &     -Q \\
%    0 &   T \end{array} \right]$.
%In particular $\Phi(D^2_{T,Q}) = QD_{t,1}$, and
%$\Phi(Q^{-k}XD^{2k}_{T,Q}) = \Phi(X)D^k_{t,1}$.
%  \end{proposition}
 \begin{proposition}
 \label{Proposition2}
    For any pair of rational nonzero parameters $(T,Q)$  and $t = T^2Q^{-1}-2$
    there is a  canonical ring isomorphism $\Phi: \mathcal R(T,Q)
    \to \mathcal R(t,1)$,
%    $\mathcal L(T,Q)$ and $\mathcal L(t,1)$,
    defined  by the conjugation $\Phi(X)=  A^{-1}XA  \in \mathcal R(t,1)$,
    where $X=[x_0\ x_1]$,
    \[
A =\left[ \begin{array}{cc} Q &     -Q \\
    0 &   T \end{array} \right]\ \
\text{and}\ \ \ T\Phi(X)=\left[ \begin{array}{cc} x_2-tQx_0 &     -Qx_0 \\
    Qx_0 &   x_2 \end{array} \right].
\]
In particular $\Phi(D^2_{T,Q}) = QD_{t,1}$, and
$\Phi(Q^{-k}XD^{2k}_{T,Q}) = \Phi(X)D^k_{t,1}$.
  \end{proposition}
  Since $t^2-4 = \Delta\ T^2/Q^2$ the respective quadratic fields
  coincide, and
 the proof of this Proposition is obtained by direct calculation.
  However its meaning can be elucidated further.
For a recursive sequence
$\{x_n\}_{n\in\mathbb Z}=X \in \mathcal R(T,Q)$
and $W = T\Phi(D_{T,Q}^{-1}) = [ -1 \ 1 ]\in \mathcal R(t)$
we have
\begin{equation}
  \begin{aligned}
  \Phi(X) =\{T^{-1}Q^{1-n}x_{2n}\}_{n\in\mathbb Z}&=Z\in \mathcal R(t), \\
  \{Q^{1-n}x_{2n-1}\}_{n\in\mathbb Z}&=T\Phi(D_{T,Q}^{-1}X) = WZ\in \mathcal R(t).
  \end{aligned}
\end{equation}
%Note that these formulas are consistent with the recursive relation
%$ Tx_{2n-1}= Qx_{2n-2}+x_{2n}$.
In particular we substitute
$X = [ 0 \ 1 ] = \{L_n\}_{n\in\mathbb Z}$ in the  formula (1), i.e.,
the identity matrix in $\mathcal R(T,Q)$, or in other words
the Lucas sequence, \cite{Wi}.
We obtain that for $I= [ 0 \ 1 ] = \{u_n\}_{n\in\mathbb Z}\in \mathcal R(t)$
and $W = \{w_{2n-1}\}_{n\in\mathbb Z}\in \mathcal R(t)$ we have
\begin{equation}
  u_n =\frac{1}{TQ^{n-1}}L_{2n}, \ \ w_{2n-1} = \frac{1}{Q^{n-1}}L_{2n-1}.
\end{equation}
(The choice of indexing for the sequence $W$ has some advantages
that will be illuminated later on.)

%For recursive sequences
%$\{x_n\}_{n\in\mathbb Z}$ and $\{y_k\}_{k\in\mathbb Z}$ of the matrices
%$D_{T,Q}$ and $D_{t,1}$ respectively, if $T\Phi([x_{0} \ x_{1}]) =
%[Qx_0\ x_2]=[y_{0} \ y_{1}]$ then $y_k =Q^{1-k}x_{2k}$. Similarly
%for the sequence $\{z_k\}_{k\in\mathbb Z}$ of the matrix $D_{t,1}$
%such that $T\Phi([x_{-1} \ x_{0}]) =
%[Qx_{-1}\ x_1]=[z_{0} \ z_{1}])$ we have  $z_k =Q^{1-k}x_{2k-1}$.

Hence one recursive sequence $\{x_n\}_{n\in\mathbb Z}$ for the parameter pair
$(T,Q)$ is split into two recursive sequences
$\{Q^{1-k}x_{2k}\}_{k\in\mathbb Z}$ and $\{Q^{1-k}x_{2k-1}\}_{k\in\mathbb Z}$
for the parameters $(t,1)$, essentially the even and odd numbered elements
of the original sequence.

We say that  parameter pairs $(T,Q)$ and $(\widetilde T, \widetilde Q)$
are {\it similar}
 if $T^2\widetilde Q = \widetilde T^2Q$.
Clearly the similarity
is an equivalence relation.
 In particular
$(T,Q)$ and $(aT,a^2Q)$ are similar  for any rational  $a$.
Note that  $(T,Q)$ and $(-T,Q)$ are always
similar, while $(T,Q)$ and $(T,-Q)$ are never  similar
(we exclude the values $T=0$ or $Q=0$ from consideration).

 Similar parameter pairs give the same
 parameter $t =T^2/Q -2$, and hence in view of Proposition
 \ref{Proposition2} they
 produce essentially the same recursive sequences. More precisely
 if  $\{x_n \}_{n\in\mathbb Z}$ is a recursive sequence
 for the parameters $(T,Q)$
 then the sequence $\{\widetilde x_n=a^{n}x_n\}_{n\in\mathbb Z}$
is a recursive sequence for the parameters
 $(\widetilde T,\widetilde Q)=(aT,a^2Q)$.
In particular
   the sequence
    $\{(-1)^nx_n\}_{n\in \mathbb Z}$ is
    a recursive sequence for the parameters   $(-T,Q)$.

We say that a parameter pair $(T,Q)$  is
{\it simple} if $T$ and $Q$ are integers and for any prime divisor
$p$ of  $T$ the parameter
$Q$ is not divisible by $p^2$. Every $(T,Q)$ is similar to  exactly two
simple pairs of the form $(\pm aP,aR)$ with  co-prime integers
$R$ and $P$,
and $R$ and $a$.  It is obtained directly by simplifying the fraction
$\frac{T^2}{Q}= \frac{aP^2}{R}$, with a square-free  $a$.
In the rest of the paper we will consider either simple pairs
$(T,Q)=( aP,aR)$ or one rational parameter $t$.

The idea to allow the rational parameter $t$,
  and achieve determinant $1$, appeared first in \cite{Wo},
  and it lead to substantial simplifications. In this paper
  we find further applications of this method.
  For any rational $t$ we introduce the notation $D_t= D_{t,1}$,
  and
  $\mathcal R(t) =\mathcal R(t,1)$.
  We will study mainly the
sequences in $\mathcal R(t), \ t=\frac{aP^2}{R}-2$. The information
obtained can be then translated into the original space $\mathcal R(T,Q)$.

%and sequences with integer initial values $[x_0 \ x_1]$.

%  This isomorphism
%  amounts to changing a recursive sequence
%$\{w_n\}_{n\in\mathbb Z}$ of the matrix $D_{T,Q}$ into the
%  recursive sequence $\{w_{2k}Q^{-k}\}_{k\in\mathbb Z}$ of the
%  matrix $D_{t,1}$.

  There is a hidden symmetry of the problem, which is
  revealed by the passage to
  the parameter $t$. The following proposition is established by inspection.
  \begin{proposition}
 \label{Proposition3}
    For any rational  $t$ there is a  canonical ring isomorphism
    $\Psi: \mathcal R(t) \to \mathcal R(-t)$
    given by the transposition
    $\Psi(X) =  X^T$.

    In particular $D_{-t}= -D_t^T$, and for a sequence
    $\{x_n\}_{n\in \mathbb Z} = X\in\mathcal R(t)$, the sequence
    $\{(-1)^{n-1}x_n\}_{n\in \mathbb Z}=\Psi(X)\in\mathcal R(-t)$.
  \end{proposition}
We say that the parameters $t$ and $-t$ are {\it twins},
  and that two parameter pairs $(T,Q)$ and
  $(\widehat T,\widehat Q)$ are {\it twins} if the respective parameters
$t$ and $\widehat t$ are twins. It is straightforward
  that two simple parameter pairs are twins if and only if they are equal to
  $(aP,aR)$ and   $(bS,bR)$, respectively,  satisfying $aP^2+bS^2 =4R$.
  Every simple parameter pair has a simple twin, essentially unique.
  For example the twin of the Fibonacci pair $(1,-1)$ is $(5,5)$.
  In general with $\Delta = T^2-4Q$ the parameter pairs
  $(T,Q)$ and $(\Delta, -\Delta Q)$  are twins
  (here we do not assume or claim their simplicity).

  Proposition \ref{Proposition2} can be generalized to ``higher powers''.
  We recall that
a recursive sequence $\{x_n\}_{n\in\mathbb Z}$ of $D_t$ can
be written in terms of the Chebyshev polynomials of the second  kind
$U_n(t)$, \cite{Wo},
\[
D^n_t= \left[ \begin{array}{cc} -U_{n-1}(t) &     -U_{n}(t) \\
    U_{n}(t) &   U_{n+1}(t) \end{array} \right], \ \ \
x_n =U_n(t) x_1-U_{n-1}(t)x_0.
\]
In particular the identity matrix $I$ of the ring $\mathcal R(t)$ represents
the sequence $\{U_n(t)\}_{n\in\mathbb Z}$.
The Chebyshev polynomials of the first kind are equal to
$C_n(t) = tr \ D^n_t$.  The sequence
$\{C_n(t)\}_{n\in\mathbb Z}$ is equal to $C =[2 \ t ]$.

It is  useful to introduce also the Chebyshev polynomials
of the third and fourth kind, \cite{Y}, namely
$W_{2k-1} = U_k+U_{k-1}$ and $V_{2k-1} = U_k-U_{k-1}$.
The sequence $\{W_{2k-1}(t)\}_{n\in\mathbb Z}$ is equal to $W = [ -1 \ 1 ]$
and the sequence $\{V_{2k-1}(t)\}_{n\in\mathbb Z}$ is equal to $V = [ 1 \ 1 ]$.
The special indexing of these polynomials simplifies the formulation
of their properties, cf. \cite{Wo}.

It is not hard to prove that
  the only rational zeroes of the polynomials $U_n$ are
  $t=0,\pm 1,\pm 2$, \cite{Wo}. We will exclude these values of the
  parameter $t$ from further considerations.
%  \begin{proposition}
%   For any natural $r$ there is a  canonical ring isomorphism
% $\Phi_r: \mathcal R(t)
%\to \mathcal R(t_r)$,
% where $t_r = C_r(t)$,  given by the conjugation
%
% $\Phi_r(X) =A^{-1}XA \in \mathcal R(t_r)$,
%    where
%$
%A =\left[ \begin{array}{cc} 1 &     -U_{r-1}(t) \\
%    0 &   U_{r}(t) \end{array} \right]$.
%In particular $\Phi_r\left(D^r_t\right) = D_{t_r}$ and
% $\Phi_r\left(XD^{kr}_t\right) = \Phi_r(X)D^k_{t_r}$.
%In particular $A^{-1}D^r_tA = D_{t_r}$.
%  \end{proposition}
  \begin{proposition}
   \label{Proposition4}
   For any natural $r$ there is a  canonical ring isomorphism
 $\Phi_r: \mathcal R(t)
\to \mathcal R(t_r)$,
 where $t_r = C_r(t)$,  given by the conjugation
 $\Phi_r(X) =A^{-1}XA \in \mathcal R(t_r)$,
    where $X=[x_0\ x_1]$,
\[
A =\left[ \begin{array}{cc} 1 &     -U_{r-1}(t) \\
    0 &   U_{r}(t) \end{array} \right]\ \
\text{and} \ \ U_r(t)\Phi_r(X) =
\left[ \begin{array}{cc} x_r-t_rx_0 &     -x_{0} \\
   x_0  &  x_{r}  \end{array} \right].
\]
In particular $\Phi_r\left(D^r_t\right) = D_{t_r}$ and
$\Phi_r\left(XD^{kr}_t\right) = \Phi_r(X)D^k_{t_r}$.
%In particular $A^{-1}D^r_tA = D_{t_r}$.
  \end{proposition}
  \begin{proof}
    For a recursive sequence
    $\{x_n\}_{n\in\mathbb Z}$ of the matrix $D_t$
    the sequence
$\{x_{kr}\}_{k\in\mathbb Z}$ is a recursive  sequence of the
  matrix $D_{t_r}$. Indeed we have
\[
    [x_0 \  x_{r}]=[x_0 \   x_{1}]
    \left[ \begin{array}{cc} 1 &     -U_{r-1}(t) \\
        0 &   U_{r}(t) \end{array} \right], \ \
         [x_r \  x_{2r}]=[x_0 \  x_{r}]A^{-1}D^r_tA,
         \]
         It follows that $A^{-1}D^r_tA = D_{t_r}$.

         Clearly  $\Phi_r$ is a ring monomorphism
         into $\mathcal R(t_r)$. Since for $X= [x_0\ x_1]$
         we have
         $U_r(t)\Phi_r(X) = [x_0\ x_r]$ and for a given $x_0$ the
         element $x_r$ assumes arbitrary rational value, then
         $\Phi$ is onto $\mathcal R(t_r)$.
         \end{proof}
  It follows from this Proposition that
  splitting a recursive sequence
  $\{x_n\}_{n\in\mathbb Z}$ of the parameter $t$  into $r$
  subsequences $\{x_{kr+j}\}_{k\in\mathbb Z}, j=0,1, \dots, r-1$,
  we obtain recursive sequences of the parameter $t_r$.

  We say that a parameter $t$ is { \it r-primitive} if there is no rational
$u$, such that $t = C_r(u)$, and it is { \it primitive} if
it is { \it r-primitive} for all natural $r \geq 2$.

%In particular  $t$ is 2-primitive if and
%only if $2+t$ and $2-t$ are not rational squares.

%Rankin,\cite{R}, studied
%{\it generators}, which are essentially   primitive parameters $t$,
%and $2$-primitivity plays an important role  in \cite{Wo}.

In view of Proposition \ref{Proposition4} we will mostly restrict
our attention to
primitive parameters (excluding again the values $t=0,\pm 1,\pm 2$,
which are not primitive).

We look now at  matrices in $\mathcal R(t)$
with determinant $1$.

 \begin{lemma}
 \label{Lemma5}
Given a rational $t$, for a rational $b$ there is
a matrix $X \in \mathcal R(t)$
with $tr \ X =b, \ \det X =1$ if and only if   $ \frac{b^2-4}{t^2-4}$
is a rational square.

The matrix $X$, if it exists, is unique up to the automorphism
of $\mathcal R(t)$.
 \end{lemma}
 \begin{proof}
We  seek an   element $X=[x_0 \ x_1] \in \mathcal R(t)$ with $\det X =1$
  and $tr \ X = b$. By a direct calculation we arrive at
  the equivalent conditions of the form
  $x_0^2= \frac{b^2-4}{t^2-4}$ and $2x_1 = tx_0 +b$.

  When there is such a rational $x_0$ then $-x_0$ gives us another matrix
  $X^{-1}$ with the same trace $b$, and there are no other in $\mathcal R(t)$.
 \end{proof}
 Clearly the condition in the Lemma  is equivalent to
  the rings $\mathcal R(t)$ and $\mathcal R(b)$ being isomorphic.

\section{Recombination of sequences}
We have seen that a recursive sequence can be split into two, or more,
recursive sequences for a new parameter.
In this section we will inspect the opposite, how to recombine
two, or more sequences of one parameter, into  recursive sequences
of another parameter.

For a simple parameter pair $(T,Q)$ and $t=T^2/Q-2$ let us recall
the special sequences in $\mathcal R(t)$ we have introduced earlier
\[
C= [2 \ t] = \frac{T}{Q}\Phi\left([2 \ T]\right),
\   W = [-1\ 1] = T\Phi\left( D^{-1}_{T,Q}\right),
\  \ V = [1\ 1] = \frac{1}{t+2} WC.
\]
%For a sequence $X\in \mathcal R(t)$ we say that
%$X$ and $CX$ are {\it polar}
%sequences, \cite{L1}.
In $\mathcal R(t)$ the Proposition \ref{Proposition6} becomes
\begin{proposition}
 \label{Proposition61}
  For any  $[z_0 \ z_1] = Z \in \mathcal R(t)$ we have
  \[
  Z^2 =
  \left[ \begin{array}{cc}
      z_0(z_1-z_{-1}) &   (z_1+z_0)(z_1-z_0) \end{array} \right],
  \]
  The even numbered elements of $Z^2$ factor into the elements
  of $Z$ and
  $CZ  = [z_{1}-z_{-1} \ z_2-z_0]$. The odd numbered elements of
  $Z^2$ factor into the elements of
  $WZ = [z_0+z_{-1} \ z_1+z_0]$ and   $VZ = [z_0-z_{-1} \ z_1-z_0]$.
  \end{proposition}

It follows from the  formula (1)
that the  sequences $Z \in  \mathcal R(t)$ and $WZ$  can be
recombined   to give a  sequence
in $\mathcal R(T,Q)$ with even-numbered elements essentially
equal to the elements of $Z$, and odd-numbered elements essentially
equal to the elements of $WZ$.

It turns out that polar duality of sequences in $\mathcal R(T,Q)$
is connected to twin symmetry in parameters. This is the subject
of the next Proposition. The twin symmetry in one parameter $t$
is described by the isomorphism $\Psi: \mathcal R(t) \to \mathcal R(-t)$ of
Proposition \ref{Proposition3}. It gives us
\[
\Psi(C_t) = -C_{-t},\   \Psi(W_t) = V_{-t}, \  \Psi(V_t) = W_{-t}.
\]
It follows that while   sequences $Z$ and $ WZ$ from $\mathcal R(t)$
can be recombined
into one  sequence $X\in \mathcal R(T,Q)$, the sequences
$Z$ and $ VZ$  can be recombined
into one  sequence $\widehat X\in \mathcal R(\widehat T, \widehat Q)$,
where
$(\widehat T,\widehat Q)= (\Delta, -\Delta Q)$ is the twin  pair
(in particular $\widehat T^2/\widehat Q = 2-t$).
More precisely we have the following Proposition.
\begin{proposition}
 \label{Proposition7}
  For  any  sequence $\{z_n\}_{n\in\mathbb Z}=Z\in\mathcal R(t)$,
  the polar sequences $\{x_n\}_{n\in\mathbb Z}=X=\Phi^{-1}(Z)\in\mathcal R(T,Q)$ and
  $\{y_n\}_{n\in\mathbb Z}=CX\in \mathcal R(T,Q)$ are given by
 \[
 \begin{aligned}
x_{2k} = TQ^{k-1}z_k,
\ \ \
&x_{2k-1} =  Q^{k-1}\left(z_k+z_{k-1}\right),\\
   y_{2k} = Q^{k}\left(z_{k+1}-z_{k-1}\right),
\ \ \
&y_{2k-1} =  TQ^{k-1}\left(z_k-z_{k-1}\right), \ k= 0,\pm 1,\dots.
\end{aligned}
\]
Further
the sequences $\{\widehat x_n\}_{n \in \mathbb Z}$ and
$\{\widehat y_n\}_{n \in \mathbb Z}$ given by
 \[
 \begin{aligned}
\widehat x_{2k} = T^{-1}\Delta^{k}x_{2k},
\ \
&\widehat x_{2k-1} = T^{-1}\Delta^{k-1} y_{2k-1}, \\
\widehat y_{2k} =\Delta^{k} y_{2k},
\ \
&\widehat y_{2k-1} =\Delta^{k} x_{2k-1}, \ \ \ \ \ \ \ k= 0,\pm 1,\dots,
\end{aligned}
 \]
% \[
% \begin{aligned}
%&\widehat x_{2k} = (-1)^k\widehat T\widehat Q^{k-1}z_k,
%\ \
%\widehat x_{2k-1} = (-1)^k \widehat Q^{k-1}
%\left( z_k- z_{k-1}\right), \\
%&\widehat y_{2k} =
%(-1)^{k+1}\widehat Q^{k-1}\left(z_{k+1}-z_{k-1}\right),
%\
%\widehat y_{2k-1} =
%(-1)^k\widehat T\widehat Q^{k-2}\left(z_k+z_{k-1}\right),
%\end{aligned}
%\]
 are polar sequences in  $\mathcal R(\widehat T, \widehat Q)$.
\end{proposition}
\begin{proof}
  The form of the even-numbered elements of $X$ follows directly from
  the formula (1).
  Once these are established we  get
  the odd-numbered elements from the recursion relation
  $Tx_{2k-1} = x_{2k}+Qx_{2k-2},\  k = 0,\pm 1, \dots ,$
  for any sequence  $\{x_n\}_{n\in\mathbb Z}\in\mathcal R(T,Q)$.

  To get the elements of $CX =\{y_n\}_{n\in\mathbb Z}$ we use
  the connection between polar sequences
  $y_n = x_{n+1}-Qx_{n-1}, n= 0,\pm 1, \dots $ and substitute
  into it the obtained formulas for $x_n$. That gives us the first part
  of the claim.

  Further we consider the sequence
  $\{(-1)^{n-1}z_n\}_{n\in\mathbb Z}=\Psi(Z)\in\mathcal R(-t)$
  and apply to it the first part. We get
  \[
  \begin{aligned}
  &\widehat x_{2k} = \widehat T \widehat Q^{k-1}(-1)^{k-1}z_k =
  \Delta^{k}   Q^{k-1}z_k =
  T^{-1}  \Delta^{k}  x_{2k},\\
  & \Delta \widehat x_{2k-1}=
  \widehat T\widehat x_{2k-1} = \widehat x_{2k} +\widehat Q\widehat x_{2k-2}
  = T^{-1}\Delta^k\left(x_{2k} -Qx_{2k-2}\right) =
  T^{-1}  \Delta^{k} y_{2k-1}.
  \end{aligned}
  \]
  We finish using the connection between polar sequences in
 $\mathcal R(\widehat T, \widehat Q)$.
 \[
  \begin{aligned}
    &\widehat y_{2k} = \widehat x_{2k+1} -\widehat Q \widehat x_{2k-1}
    = T^{-1}  \Delta^{k}  \left(y_{2k+1} +Qy_{2k-1}\right) =
  \Delta^{k}   y_{2k},\\
  &\widehat y_{2k-1} = \widehat x_{2k} -\widehat Q \widehat x_{2k-2}
=T^{-1}  \Delta^{k}  \left(x_{2k} +Qx_{2k-2}\right) =
  \Delta^{k} x_{2k-1}.
  \end{aligned}
  \]
\end{proof}
The pairs of sequences $Z,WZ $
and $CZ, VZ\in \mathcal R(t)$ correspond to
two polar sequences   in $\mathcal R(T,Q)$.
It transpires from Proposition \ref{Proposition7} that
essentially exchanging odd numbered elements in polar sequences
in $\mathcal R(T,Q)$ gives us polar sequences in
$\mathcal R(\widehat T, \widehat Q)$.
Note that since $\Phi(C_{T,Q}) = T^{-1}QC_t$, the polarity of two
sequences is preserved under the isomorphism
$\Phi:\mathcal R(T,Q)\to \mathcal R(t)$ only up to
scalar factors.

Guided by this recombination phenomenon we want to consider
two rational parameter values $t\neq \pm a$ such that
$C_r(t) = \pm C_r(a)$, for a prime  $r$.
It turns out that it is
fairly rare.
 \begin{lemma}
 \label{Lemma8}
   The equation  $C_r(t) = \pm C_r(a)$, for a prime  $r$,
   has a rational solution $t\neq \pm a$ if and only if $r=2$ or $3$.

   The equation $C_2(t) = - C_2(a)$
   is equivalent to $t^2+a^2 = 4$.

   The equation  $C_3(t) =  C_3(a)$ is solvable if and only if
   $t^2-4 = -3f^2$, for some rational $f$, with $a = \frac{-t\pm 3f}{2}$.
 \end{lemma}
 \begin{proof}
   It is a direct check that the solutions listed are indeed solutions.
In particular we have $C_3(t) -C_3(a)=(t-a)\left(t^2+ta +a^2-3\right)$,
and solving the quadratic equation gives us the answer above.

Let us now assume that
there is such a rational solution $t\neq a$  for a  prime $r$.
We substitute $x=t$ and $x=a$ into the following identity valid
for Chebyshev polynomials
for any natural  $n$, to obtain
\[
C_n^2(x) -(x^2-4)U_n^2(x) =4, \ \ \frac{a^2-4}{t^2-4} =\frac{U_r^2(t)}{U_r^2(a)}.
\]
Hence by Lemma \ref{Lemma5} we get an   element
$X=[x_0 \ x_1] \in \mathcal R(t)$
with $\det X =1$ and $tr \ X = a$.
We have that $X^r$ and $D_t^r$
have the same traces, and determinant $1$.
It follows by Lemma \ref{Lemma5} that $X^{\pm r}=  D_t^{r}$,
and further  $\left(X^{\mp 1}D_t\right)^r= I$.

It is well known that there are roots of unity of  order $n$ in
$SL(2, \mathbb Q)$ only for $n=2,3,4$ and $n=6$, which ends
the proof. (We will give an independent proof of a more general fact
in Proposition \ref{Proposition17}.)
 \end{proof}
 Note that the equation $C_2(t) = - C_2(a)$ is effectively equivalent
 to $C_4(t) = C_4(a)$.

 We will refer to the $r=2$ case as {\it circular},
 and $r=3$ as  {\it cubic}, jointly as {\it cyclotomic} cases.
 Both cases appear explicitly
 in Ballot, \cite{B}, in the language of two parameters.
 For a given parameter $t$, the other parameters $a$, if present,
will be called the {\it associate} parameters.
Note that the circular case and the cubic case are
exclusive; a parameter $t$  cannot be circular and cubic
at the same time. By the above identity for Chebyshev polynomials
we get that for any natural $n$ and rational $t$ the parameter
$C_n(t)$ is cyclotomic
if and only if $t$ is cyclotomic, with the preservation of the
circular and cubic type.

 For associate parameters $t$ and $a$ there is a canonical ring isomorphism
 $\Theta=\Theta_{t,a}:\mathcal R(t) \to \mathcal R(a)$, which is defined
 on the basis of Propositions \ref{Proposition3} and \ref{Proposition4}.
 So
 $\Theta_{t,a} =\Phi^{-1}_{2,a}\circ \Psi \circ \Phi_{2,t}$ in the circular case,
 and
 $\Theta_{t,a} =\Phi^{-1}_{3,a}\circ  \Phi_{3,t}$ in the cubic case.

 The next Proposition details the relation between the
 sequences  $X$ and $\Theta(X)$
 in the  circular case and the recombination of the polar sequences
 $X$ and $CX$ from $\mathcal R(t)$ into polar sequences in
  $\mathcal R(a)$.
 \begin{proposition}
  \label{Proposition9}
   In the circular case for the associate parameters $t$ and $a$,
   for any sequence $\{x_n\}_{n\in\mathbb Z} =X \in \mathcal R(t)$
   we have
   \[
   \Theta(X) = t^{-1}[-ax_0 \ tx_1-2x_0] = t^{-1}[-ax_0 \ x_2-x_0].
   \]
   In particular
   \[
   \Theta(C_t) = -at^{-1}C_a, \  \Theta(D_tC_t) = -aD_a, \
   \Theta(D^2_t) = -D^2_a.
   \]
   Further, if $\{y_n\}_{n\in\mathbb Z}=CX$,
   $\{\widehat x_n\}_{n\in\mathbb Z} = \Theta(X)$
   and    $\{\widehat y_n\}_{n\in\mathbb Z} = C_a\Theta(X)$
   then
   \[
   \begin{aligned}
&\widehat x_{2k} = (-1)^{k-1}at^{-1}x_{2k},
\ \
\widehat x_{2k+1} =(-1)^kt^{-1}y_{2k+1},\\
&\widehat y_{2k} = (-1)^{k}y_{2k},
\ \
\widehat y_{2k+1} =(-1)^kax_{2k+1},\ \ \ \
k= 0,\pm 1,\dots.
\end{aligned}
\]
 \end{proposition}
The proof is by a straightforward calculation.
Note the parallel of the last Proposition with Proposition \ref{Proposition7}.
Again we learn that essentially exchanging odd numbered elements in a polar
pair in $\mathcal R(t)$ we obtain a polar pair in $\mathcal R(a)$.

The  recombination phenomenon is more complicated in the cubic case
in which $\mathcal R(t)$ is isomorphic to $\mathbb Q(\sqrt{-3})$.
To describe  it we employ the two roots of unity of order $3$,
$S =S_t =-\frac{1}{2f} [2\  t+f]$ and $R=R_t =
\frac{1}{2f}[2\ t-f]$ in $\mathcal R(t)$.
We have $\det S =1=  \det R, tr \ S = -1 = tr \ R, S^3 =  I = R^3$
and $SR =  I$.
\begin{proposition}
 \label{Proposition11}
   In the cubic  case for the associate parameters $t$ and $a$,
   for any sequence $\{x_n\}_{n\in\mathbb Z}=X \in \mathcal R(t)$
   we have
   \[
  \left(t^2-1\right) \Theta(X) = [(a^2-1)x_0 \ (a-t)x_0+(t^2-1)x_1 ]
   = [(a^2-1)x_0 \  ax_0+x_3].
   \]
   In particular
   \[
   \Theta(D_tS_t) = D_a, \    \Theta(D^{-1}_tR_t) = D^{-1}_a.
   \]
   Further if $\{y_n\}_{n\in\mathbb Z}=SX,\ \{z_n\}_{n\in\mathbb Z}=RX$ and
   $\{\widehat x_{n}\}_{n\in\mathbb Z}=\Theta(X)$  then
 \[
\widehat x_{3k} = cx_{3k},\
 \widehat x_{3k+1} = cy_{3k+1}, \
 \widehat x_{3k-1} = cz_{3k-1}, \ \ c = \frac{a^2-1}{t^2-1},\
 k= 0,\pm 1,\dots.
 \]
\end{proposition}

\begin{proof}
  The element $D_t^3 \in \mathcal R(t)$ has exactly three
  roots   of order $3$, namely  $D_t, D_tS, D_tR$.
  Since $\Phi_{3,t}(D_t^3) = D_{t_3} = D_{a_3} = \Phi_{3,a}(D_a^3)$,
  where $t_3= C_3(t) = C_3(a) = a_3$, we conclude
  that
  \[
  \{\Theta(D_t),\Theta (D_tS_t), \Theta(D_tR_t)\}=
  \{D_a, D_aS_a, D_aR_a\}.
  \]
  Since all of these matrices have determinant $1$, and
  $\Theta$ preserves  traces, it remains to compare $tr \ D_tS_t$ and
  $tr \ D_a=a$.
  We have $tr \ D_tS_t = \frac{4-t^2-tf}{2f} =\frac{3f^2-tf}{2f} =a $.
  Hence $\Theta(D_tS_t) = D_a$. It follows that
  $\Theta(D^{-1}_tR_t) =\Theta\left((D_tS_t)^{-1}\right)= D^{-1}_a$,
  since $\Theta$ is a ring isomorphism.

  The first equality  is obtained  by direct calculation.
  Taking into account that  $\Theta(D_t^3)= D_a^3 $ and
  $\Theta(D_t^3X)= D_a^3\Theta(X)$ we obtain
  $\widehat x_{3k} = c x_{3k}$.

To get $\widehat x_{3k+1} = c y_{3k+1}$ we observe that
\[
\Theta(X)D_a^{3k+1} = \Theta(X) \Theta(D_tS_t) \Theta(D_t^{3k}) =
\Theta(S_tXD_t^{3k+1}).
\]
In a similar fashion we arrive at $\widehat x_{3k-1} = c z_{3k-1}$.
\end{proof}
A proof of  the last equality in the Proposition
by direct calculation requires the application of the following
surprising identity
$\frac{a^2-1}{t^2-1} = -\frac{2f}{t+f}$.

Note that for any recursive sequence   $\{x_n\}_{n\in\mathbb Z}$,
in $\mathcal R(t)$, we have

 $ (t^2-1)x_{3k+1} =  \left(tx_{3k}+x_{3k+3}\right)$,
 $(t^2-1)x_{3k-1}= \left(x_{3k-3}+tx_{3k}\right), \ k= 0,\pm 1,\dots
$. Using these identities we can reformulate the
last Proposition somewhat differently.
\begin{proposition}
 \label{Proposition12}
  In the cubic case for the associate parameters $t$ and $a$,
  for any sequence $\{x_n\}_{n\in\mathbb Z}=X \in \mathcal R(t)$
the sequence   $\{\widehat x_n\}_{n\in\mathbb Z}$ given by
\[
\widehat x_{3k} = x_{3k}, \ \widehat x_{3k+1} =  \frac{1}{a^2-1}\left(ax_{3k}+x_{3k+3}\right),
 \ \widehat x_{3k-1}= \frac{1}{a^2-1}\left(x_{3k-3}+ax_{3k}\right),
 \]
 $k= 0,\pm 1,\dots$ belongs to $\mathcal R(a)$.
 Further if $ \{y_n\}_{n\in\mathbb Z}=SX$ and  $\{z_n\}_{n\in\mathbb Z}=RX$ then
 \[
 \widehat x_{3k+1} =  y_{3k+1}, \  \widehat x_{3k-1} = z_{3k-1}, \ \
 k= 0,\pm 1,\dots.
 \]
\end{proposition}

Hence in the cubic case $X$ and $\Theta(X)$ essentially
share the same
 elements with indices divisible by $3$.
Further, essentially exchanging appropriate subsequences
 in the sequences $X, SX$ and $RX$ in $\mathcal R(t)$ we obtain
 three recursive sequences $\Theta(X), \Theta(SX)$ and
 $\Theta(RX)$ in $\mathcal R(a)$.

 \section{Sequence  groups}

The following constructions depend on
the choice of $T,Q\neq 0$, which are fixed. We will not show this
dependence explicitly in notation unless it may cause ambiguity.

We are going to identify recursive sequences which differ by a scalar
rational factor. To that end
  we projectivize the ring $\mathcal R$, i.e., we consider two elements of
  $\mathcal R\setminus \{0\}$,
  $[x_0 \   x_1]$ and $[y_0 \   y_1]$, equivalent
  if there are integers $k\neq 0,l\neq 0$, such that
  $k[x_0  \  x_1]= l[y_0 \   y_1]$.
  Clearly it is an equivalence relation and we consider the space of
  equivalence classes, with the exclusion of elements $X$ with
  $\det \ X = 0$. We denote the  resulting space by $\mathcal L$.
  \begin{lemma}
   \label{Lemma13}
    The multiplication of matrices is well defined for the equivalence
    classes in $\mathcal L$, and $\mathcal L$ becomes a commutative group.
      \end{lemma}
  We will call the group $\mathcal L=\mathcal L(T,Q)$ the
  {\it sequence group}.

Let us note that every equivalence class in $\mathcal L$ contains
a matrix $X=[x_0 \ x_1]$ with  integer, relatively prime   $x_0$ and $x_1$.
We call such a matrix a {\it reduced} representative.
Clearly every equivalence class
  has exactly two reduced representatives, differing by
  the sign.

  We address now the problem of taking square roots in the sequence group
  $\mathcal L$.
\begin{proposition}
 \label{Proposition14}
  An element $[y_0 \ y_1]= Y \in \mathcal L(T,Q)$ is a square if and only if
  $\det Y$ is a square, $\det Y=\lambda^2, \lambda  \in \mathbb Q$.
  The solutions $[x_0 \ x_1]= X  \in \mathcal L(T,Q)$
  of the equation $X^2 =Y$
  are the row eigenvectors of  the matrix
  $ \left[ \begin{array}{cc} -y_1 & Qy_0-Ty_1      \\
         y_0 &   y_1 \end{array} \right]$ with the eigenvalues $\pm \lambda$.
  Moreover if $y_0\neq 0$ then the square roots are equal to
  $[y_0 \ y_1\pm \lambda]$.
\end{proposition}
\begin{proof}
  We are going to solve the matrix equation $b X^2 = Y$,
    for a rational $b$.
    We rewrite the equation as $ b X = YX^{-1}$.
    Changing the free rational parameter $b$ to another one $\lambda$
    we get
    \[
 \lambda   \left[ \begin{array}{cc} -Qx_{-1} & -Qx_0      \\
       x_0 &   x_1 \end{array} \right] =
   \left[ \begin{array}{cc} y_1-Ty_0 & -Qy_0      \\
       y_0 &   y_1 \end{array} \right]
   \left[ \begin{array}{cc} x_{1} & Qx_0      \\
       -x_0 &   x_1-Tx_0 \end{array} \right].
   \]
   Clearly if the problem has a solution
   then $\det Y =\lambda^2$.
   Further, since this is an equation in the sequence group $\mathcal L(T,Q)$,
   it is equivalent to the equality of the second rows of the matrices
   on the left, and on the right. That gives us the eigenvalue problem
   \[
     [x_0 \  x_{1}]\left[ \begin{array}{cc} -y_1 & Ty_0-Qy_1      \\
         y_0 &   y_1 \end{array} \right]
  =\lambda [x_0 \  x_{1}]
\]
The characteristic equation is $\lambda^2 =\det Y$, so that the rational
square roots
of $\det Y$ are the eigenvalues.
\end{proof}
This proposition will facilitate our calculations. It also has an interesting
corollary.
  \begin{corollary}
   \label{Corollary15}
    For the group homomorphism
    $\Upsilon : \mathcal L \to \mathcal L, \ \Upsilon(X) =
    X^2$
    the image $\Upsilon(\mathcal L)$
    contains exactly those  $Y\in\mathcal L$ for which
    there is a matrix $A \in  SL(2,\mathbb Q)$ and
      integers $k \neq 0, l \neq 0$ such that
      $kY= lA$.
        \end{corollary}
  \begin{proof}
    %  Despite the appearance the claim in the Corollary is actually constructive.
    
    If $Y\in \Upsilon(\mathcal L)$ then $Y = X^2$ for a matrix
    $X \in \mathcal R$, so that $\det Y =\lambda^2, \lambda = \det X$.
    The sought after element $A \in SL(2,\mathbb Q)$ is equal to
    $A= \frac{1}{\lambda}Y$.

    Conversely we start with the equation
    $kY= lA$ for a matrix $A \in  SL(2,\mathbb Q)$.
    It implies that $\det Y = \left(\frac{l}{k}\right)^2$.
    By Proposition \ref{Proposition14} it follows
    that there is an element $X\in \mathcal L$ such that
    $bX^2 = Y$ for some  $b \in \mathbb Q$. Hence $Y \in \Upsilon(\mathcal L)$
  \end{proof}

  From now on we restrict our attention to the one parameter
  sequence group $\mathcal L=\mathcal L(t)$.
  For $X = [x_0  \  x_1], Z = [x_0  \  x_{-1}]$ in $\mathcal R(t)$
  we have $XZ = -(\det X) I$,
  which gives us the  inverse elements in $\mathcal L(t)$.
  It follows that  if $X$ in $\mathcal L(t)$ represents
  the sequence $\{x_k\}_{k\in \mathbb Z}$
  then $X^{-1}$ can be
    represented by  the sequence $\{x_{-k}\}_{k\in \mathbb Z}$.

  %In particular $D\in \mathcal L$.

%  on the basis of the following
%  matrix identities.
    %    For $Y = [y_0  \  y_1], Z = [y_1  \  y_{0}]$ we have
    %$YZ = -(\det Y) D^{-1}$.
    %Also for $ZD = [y_0  \  y_{-1}]$ we have $Y(ZD) = -(\det Y) I$.
 %Hence if $Y$ is the sequence $\{y_k\}_{k\in \mathbb Z}$ then $Y^{-1}$ can be
%    represented by the sequence $\{y_{1-k}\}_{k\in \mathbb Z}$.

To study  the number-theoretic
properties of recursive sequences we introduce another kind
of sequence groups.
For a fixed  $t\in \mathbb Q, \ t\neq 0,\pm 1,\pm 2$, let $\Pi_t$
be the set of
odd primes not dividing the numerator, or denominator of $t$ and
$\delta =t^2-4$
($t$ is assumed to be a fully reduced simple fraction).
 For a prime $p \in \Pi_t$ we consider
the subgroup $\widetilde{\mathcal L_p}(t)\subset \mathcal L(t)$
of equivalence classes with the determinant of the reduced representative
 $\neq 0 \mod p$. There is the canonical homomorphism
$\mathcal M : \widetilde{\mathcal L_p}(t) \to \widetilde{GL}(2,\mathbb F_p)$,
where $\widetilde{GL}(2,\mathbb F_p) = GL(2,\mathbb F_p)/\mathbb F_p^*$.
We define
the {\it  sequence group $\mod  p $}, $\mathcal L_p(t)$, as a subgroup of
$\widetilde{GL}(2,\mathbb F_p)$,
$\mathcal L_p(t) =\mathcal  M\widetilde{\mathcal L_p}(t)$.

Let us note that one can equivalently introduce the group $\mathcal L_p(t)$
by repeating the construction of $\mathcal R(t)$  with the
replacement of recursive  sequences with rational elements
by recursive sequences with elements in $\mathbb F_p$.

%  For the ease of notation we employ the Legendre symbol $(a|p) =1$
%for a rational $a$,
%and an odd prime $p$ not dividing the denominator of $a$; it is equal to
%$1$ if $a\ \mod \ p$ is  a square residue, equal to $-1$
% if $a\ \mod \ p$ is a square non-residue,
%  and equal to $0$
%  if the numerator of $a$ is divisible by $p$.
\begin{theorem}
 \label{Theorem16}
  The        sequence group            $\mathcal L_p(t)$ is cyclic
  of order $p-1$ if $\delta = t^2-4$ is a square  $\mod  p$,
  and of order $p+1$ if $\delta$ is not a square $\mod  p$.
\end{theorem}
\begin{proof}
  For the proof we will trace back the construction of
  $\mathcal L(t)$, repeating it $\mod  p$. So we consider
  the matrix $D_t$ as an element in $SL(2,\mathbb F_p)$, and
  we obtain $\mathcal R_p$, the ring of $2\times 2$ matrices
  with
  entries from $\mathbb F_p$, commuting with $D_t$.
  If $\delta = t^2-4$ is a square  $\mod  p$ then
  the matrix $D_t$ has two different eigenvalues, and it can
  be diagonalized.
  The ring $\mathcal R_p$ is then isomorphic to the ring
  $\mathbb F_p\times\mathbb F_p$.
  The projectivization of $\mathcal R_p$ has $p+1$ elements,
  however two of them have
  zero determinants. It leads to the cyclic group $\mathcal L_p$ of order
  $p-1$.

If $\delta$ is not a square $\mod  p$ then
the matrix $D_t$ has no eigenvalues, and  it is well known that
in such a case the  ring  of matrices of the form $aI+bD_t,
a,b\in \mathbb F_p$
%(coinciding with $\mathcal R_p$)
is isomorphic to the finite field  $\mathbb F_{p^2}$.
The projectivization $\mathcal L_p$ is then
isomorphic to the quotient group $\mathbb F_{p^2}^*/\mathbb F_{p}^*$,
which is cyclic of order $p+1$.
\end{proof}
As an application of the sequence groups $\mathcal L_p$ we prove the
following useful fact.

\begin{proposition}
 \label{Proposition17}
   For any rational $t$
   there are no elements of finite order $k \geq 2$ in $\mathcal L(t)$,
   except for $k = 2,3,4,6$.
\end{proposition}
\begin{proof}
  If there is such an element $X=[x_0 \ x_1]$ then, assuming
  that $X$ is a reduced representative, for all odd primes
  $p\in \Pi_t$ not dividing $x_0$ and $\det X\neq 0 \mod p$,
  the  homomorphism
  $\mathcal M : \widetilde{\mathcal L_p}(t) \to \mathcal L_p$
  takes $X$ into an element of order $k$.
  It follows that $k$ divides the order of $\mathcal L_p$, equal to  $p\pm 1$.
  Hence $p = \pm 1  \mod  k$ for all $p$ with only finitely many exceptions.
  It is a contradiction except when $k=2,3,4,6$.
\end{proof}

We end this section with a  discussion of the torsion subgroup of
$\mathcal L(t)$.  Let us  note that   all
the isomorphisms of the rings discussed in Section 2 translate into
    isomorphisms of respective sequence groups.

Using  Proposition \ref{Proposition14} we can
easily enumerate elements of order $2^k$.
First of all we  calculate  square roots of identity in $\mathcal L(t)$
to conclude that there is exactly one  element  of order $2$,
namely   $C= [2 \ t]$. We have $C^2 = (t^2-4)I$.

   To get elements of order $4$ we need to take square roots of $C$.
   Since $\det C = 4-t^2$ it is possible only in the circular case.
\begin{proposition}
 \label{Proposition18}
  There are elements of
  order $4$ in $\mathcal L(t)$
  only in the circular case, i.e.,  for $t$ such that $t^2-4=-a^2$,
  for a rational
  $a$.  The  elements of order $4$ are then $G=[2 \ t+a]$
  and $H=[2 \ t-a]$.

For any rational $t$  there are  no elements of order $8$ in $\mathcal L(t)$.
\end{proposition}

\begin{proof}
  Elements of order $4$ must be square roots of $C$.
  Assuming $4-t^2 =a^2$, for a rational $a$, we solve the matrix
  equation $a X^2 = C$ using Proposition \ref{Proposition14}.
We obtain the only two elements of order $4$ in $\mathcal L(t)$,
$G = [2 \ t+a]$ and $H = [2 \ t-a]$. We have   $G^2 = 2a C$, $H^2= -2a C$
and $GH = -2a^2 I$.

Since $\det G = 2a^2=\det H$, these elements do not have square roots.
\end{proof}
In the  circular case the torsion subgroup of  $\mathcal L(t)$
is isomorphic to $\mathbb Z_4$.

We now address the elements of odd order $r$.
  \begin{proposition}
   \label{Proposition19}
    There are elements of odd order $r$ in $\mathcal L(t)$ only
    in the cubic case, i.e., for
    $r=3$ and $t$ such that $t^2-4=-3f^2$ for a rational
    $f$. The elements are $S=[2 \ t+f]$ and $R=[2 \ t-f]$.
    Furthermore, in this case there are also two elements of order $6$,
$Y = [2 \ t+3f]$ and $Z = [2 \ t-3f]$.
  \end{proposition}
\begin{proof}
  By Proposition \ref{Proposition17} we only need to find all elements of
  order $3$.

  Let $X=[x_0 \  x_{1}] \in \mathcal L(t)$ be an element of  order $3$.
    It means that there is $\lambda\in\mathbb Q$ such that $X^3=\lambda I$.
    It follows that $\left(\det \ X\right)^3 =\lambda^2$, and hence
    $\det X = \lambda^2\left(\det \ X\right)^{-2}$ is a rational square.
    So without loss of generality we can
    assume that $X\in SL\left(2,\mathbb Q\right)$ and $X^3= I$. The trace
    of $X$ must be a solution of the polynomial equation $C_3(u)=2$,
    and hence it is equal to $-1$. We arrive at the equations
    \[
    \det X = x_1^2-tx_1x_0+x_0^2 =1, \ \ \ tr \ X = 2x_1-tx_0=-1.
    \]
    Substituting $x_1$ in the first equation using the second equation
    we obtain
    \[
    (t^2-4)x_0^2=-3.
    \]
    If we put $t^2-4=-3f^2$ then
    $x_0 = \pm\frac{1}{f}, x_1 = \frac{\pm t-f}{2f}$.
    That gives us the only two elements of order $3$,
    $S=[2 \ t+f]$ and $R=[2 \ t-f]$.

    Elements of order $6$ must be square roots of the elements of order $3$.
Since $\det S = 4f^2 = \det R$ we obtain by Proposition \ref{Proposition14}
two elements of order 6
$Y = [2 \ t+3f]$ and $Z= [2 \ t-3f]$, with $\det Y = 12f^2 = \det Z$.
    \end{proof}
 
  We have  $S^2 = 2fR, R^2 = -2fS, SR = -4f^2 I $. Further
  $Y^2 = 6f S$, $Z^2 = -6f R$,  and $YZ = -12f^2 I$.
In the  cubic case the torsion subgroup of  $\mathcal L(t)$
is isomorphic to $\mathbb Z_6$.

%Let us note that it follows from our discussions
%that there are roots of unity of any natural order $n\geq 2$
%in some $\mathcal R(t)$ only for
%$n=2,3,4$. Indeed $-I$ has a square root only for $t=0$
%($\mathcal R(0)$ is isomorphic to $\mathbb C$)
%or in the circular case.
%Further by Proposition 17 we need only to inspect the elements
%of order $4$ and $6$. We now know all such elements in
%$\mathcal L(t)$, namely  $G,H,Y,Z$, in the cyclotomic cases, and
%it is clear that they are not of finite order in $\mathcal R(t)$.

  \section{The Laxton group and its torsion subgroup}

  Let $\mathcal D$ be the cyclic subgroup of $\mathcal L$ generated by
  the  element $D$. The Laxton group is the quotient group
  $\mathcal G = \mathcal L/\mathcal D$. The meaning of this definition
  is that while elements of $\mathcal L$ can be thought of as  recursive
  sequences, considered up to a multiplicative factor,
  in the Laxton group we also identify sequences
  that differ by a shift.

  The conjugacy of Proposition \ref{Proposition2} becomes the  homomorphism
  of $\mathcal G(t)$ onto
  $\mathcal G(T,Q)$, a $2$ to $1$ mapping.
   We denote the
  homomorphism  by $\Xi: \mathcal G(t) \to \mathcal G(T,Q)$.
However   it seems that there is no
  canonical isomorphic embedding of
  $\mathcal G(T,Q)$ into $\mathcal G(t)$.
  The explanation for this is that for
  the recursive sequence $\{x_n\}_{n\in\mathbb Z}$ of the matrix $D_{T,Q}$
  there is no
  canonical choice between the two recursive sequences
  $\{x_{2k}Q^{1-k}\}_{k\in\mathbb Z}$ and
$\{x_{2k-1}Q^{1-k}\}_{k\in\mathbb Z}$, which comprise the complete
  preimage under the homomorphism $\Xi$.

  The kernel  of the homomorphism $ \Xi$ is comprised of
  two elements, identity $I$ and $W = [-1\ 1]$. Indeed using the notation of
   Proposition \ref{Proposition2} we get
  \[
  TA^{-1}D_{T,Q}A = Q\left[ \begin{array}{cc} 1 &     -1 \\
      1 &   t+1 \end{array} \right]=
  Q\left[ \begin{array}{cc}t+ 1 &     1 \\
    -1 &   1 \end{array} \right]\left[ \begin{array}{cc} 0 &     -1 \\
      1 &   t \end{array} \right]=QWD_t.
  \]

  Similarly it follows from Proposition \ref{Proposition4} that for any
  natural $n$ we have
  the homomorphism of $\mathcal G(t_n)$ onto
  $\mathcal G(t)$, which we denote by
$\Xi_n: \mathcal G(t_n) \to \mathcal G(t)$.
  The kernel of $\Xi_n$ is obtained from the image of the set  of powers
  of $D_t$ by the ring isomorphism $\Phi_n$ of Proposition \ref{Proposition4}.

  We are going to list all elements of finite order in a Laxton group
  $\mathcal G(t)$. Clearly
  they come  from solutions
  of the matrix equation $X^r =\lambda D^k$ for  natural  $r$, integer $k$
  and rational $\lambda$. We already know all the  solutions for $k=0$.

To get all the elements of order
   $2$ in
   $\mathcal G(t)$ we need to solve the equation $ X^2 = \lambda D^{2m-1}$ for
   integer $m$.  For $Y =XD^{-m}$ we get the equation
$Y^2 = \lambda D^{-1}$.
   By Proposition \ref{Proposition14} it is solvable because $\det D =1$ and
   we obtain
   two  elements, $W=[-1 \ 1]$ and   $V=[1 \ 1]$.
We have $W^2 = (t+2)D^{-1}$ and $V^2 = (t-2)D^{-1}$.
Moreover $WC = (t+2)V$.

Further, $\det W = 2+t$ and  $\det V = 2-t$,
so that by Proposition \ref{Proposition14}
$W$ and $V$ do not have square roots if and only if both
$t$ and $-t$ are $2$-primitive. In such a case we say that
$t$ is {\it twin primitive}.

%We can thus claim that for
%twin primitive $t$ there are no elements of order $4$ in
%$\mathcal G(t)$, except in the circular case.

\begin{proposition}
 \label{Proposition20}
  For a primitive and twin primitive $t$, which is not cyclotomic,
  the torsion subgroup of the Laxton group
  $\mathcal G(t)$ contains only $4$ elements $I,C,W$ and $V$,
  and it is isomorphic to  $\mathbb Z_2\times\mathbb Z_2$.
  \end{proposition}
\begin{proof}
  We continue the search for solutions
  of
  the matrix equation $X^r =\lambda D^n$
  for an odd prime $r$, rational $\lambda$ and integer $n\neq 0$.
  We have  $(\det X)^r = \lambda^2$, and hence
  $\det X= \frac{\lambda^2}{(\det X)^{r-1}}$
  is a rational square. It is enough then to solve the equation
  $X^r=D^n$ for $X\in \mathcal R(t), \det X =1$, and $n \neq 0$.
  If $n=mr$ is divisible by $r$ then $\left(XD^{-m}\right)^r=I$,
  and the solutions are shifts by $D^{m}$ of the elements of
  odd  order $r$ in the sequence group $\mathcal L(t)$.
  By the assumption that $t$ is not cubic there are no such
  solutions.

  If $r$ does not divide  $n$   there are
  integers $k,l$ such that $ln=kr+1$. It follows that
  $X^{lr}=D^{kr+1}$, or equivalently $(X^lD^{-k})^r=D$.
  It follows that $t =C_r(u)$, where $u = tr \ X^lD^{-k}$.
  Hence contrary to the assumption $t$ is not primitive.
  We arrive at the conclusion that under the assumptions of this
  Proposition there are no elements of odd prime order in the
  Laxton group.

  To find all elements of order $4$ we need only to take square
  roots of the elements of order $2$, namely $C,W$ and $V$.
  Taking the square roots in the Laxton group $\mathcal G(t)$
  is equivalent
  by Proposition \ref{Proposition14} to taking square roots in  $\mathcal L(t)$.
  Since $t$ is twin primitive $W$ and $V$ do not have square roots.
  Further, $\det C = 4-t^2$ is  a rational square only if $t$ is circular,
  which we have excluded. Hence
  there  are  no elements of order $4$.
\end{proof}
\begin{definition}
  A cyclotomic $t$ is circular (cubic) primitive if $t$ and its associate
  values are $2$-primitive ($3$-primitive).
  \end{definition}
Let us note that for a circular $t$ and any odd prime $r$ the $r$-primitivity
is automatically shared by $t$ and its associate value.
Indeed if $t^2+a^2 =4$ then for $\widetilde C = -a^{-1}C$
we have that $\det \widetilde C =1$ and $ tr\ \widetilde CD =a$.
Hence the associate value $a$ is $r$-primitive if and only if
$\widetilde CD $ has no $r$-th root in $\mathcal L(t)$.
However, $\widetilde C^r =\pm \widetilde C$ so that $\widetilde C$
has always the  trivial  $r$-th root.

Similarly for a cubic $t$ the $r$-primitivity is shared by $t$ and
its associate values for any prime $r$ with the exception of $r=3$.
In this case the associate values of $t$ are the traces of the matrices
$SD$ and $RD$ of determinant $1$. The matrices $S$ and $R$ have trivial
$r$-th roots in $\mathcal L(t)$ except for $r =3$. Indeed, for example
the $5$-th root of $S$ is $R$, and the $5$-th root of $R$ is $S$.

\begin{proposition}
 \label{Proposition22}
  In the circular case if $t$ is primitive and circular primitive
  then the torsion subgroup
  of the Laxton group $\mathcal G(t)$
is isomorphic to $\mathbb Z_2\times \mathbb Z_4$.
  The elements $I,C,G,H$ form a subgroup    isomorphic to $\mathbb Z_4$.
  The other four elements are obtained by
the   translation by $W$.
  \end{proposition}
\begin{proof}
 % The proof is parallel to the proof of Proposition \ref{Proposition21}.
  The torsion
  subgroup of $\mathcal L(t)$ gives us the subgroup
  $\mathbb Z_4$  of $\mathcal G(t)$. The translations by $W$
  produce the subgroup isomorphic to $\mathbb Z_2\times \mathbb Z_4$.

  By the argument in the proof of Proposition \ref{Proposition20}
  there are no elements
  of odd prime order in $\mathcal G(t)$. It remains to check if
  the elements $WG$ and $WH$ of order $4$ can have square roots.
  We have  $\det WG = 2(2+t)a^2 = \det WH$, where $a$ is the associate value.
  Since for $t^2+a^2=4$ we have
$2(2+t)(2+a) = (2+t+a)^2$
  we conclude that $2(2+t)$ is a rational square if and only if
  $a$ is  $2$-primitive.
  We conclude by Proposition \ref{Proposition14} that under the assumption
  of circular primitivity there are no elements of order $8$.
\end{proof}
Let us note that we have established in the proof a simple criterion
of circular primitivity, namely a circular $t$ is circular primitive if and only
if $2+t$ and $2(2+t)$ are not rational squares.

%We say that a primitive circular $t$
%is {\it circular primitive} if $2(2+t)$
%is not a rational square. Circular primitivity is equivalent
%to $2$-primitivity of the associate parameter $a$.
%Indeed we have
%\[
%2(2+t)(2+a) = (2+t+a)^2.
%\]
%Let us note that for a circular $t$ for any odd prime $r$ the $r$-primitivity
%of $t$ and $a$ are equivalent. Indeed for $u,v$ such that $u^2 +v^2 =4$
%we have for any odd $n$ the following identity for Chebyshev polynomials,
%cf. \cite{Wo},
%\[
%C_n^2(u)+C_n^2(v) = 4.
%\]

\begin{proposition}
 \label{Proposition21}
  In the cubic case if $t$ is primitive and cubic primitive then
  the torsion subgroup of the Laxton group
  $\mathcal G(t)$ is isomorphic to $\mathbb Z_2\times \mathbb Z_6$.
  The elements $I,C,S,R,Y,Z$ form a subgroup
  isomorphic to $\mathbb Z_6$. The other six elements
  are obtained by the   translation by $W$
\end{proposition}
\begin{proof}
  The $6$
  elements comprising the torsion subgroup of $\mathcal L(t)$
  from Proposition  \ref{Proposition19},  give us also a subgroup,
  when treated as elements of $\mathcal G(t)$.

  The argument in the proof of Proposition  \ref{Proposition20} applies,
  and we can
  conclude that there are no elements of odd prime order $r$ in $\mathcal G(t)$,
  which do not come from finite order elements in $\mathcal L(t)$.

  We still need to look for roots in $\mathcal G(t)$
  of the elements of the torsion subgroup of $\mathcal L(t)$.
  Since cubic case excludes the circular case  there are
  no elements of order $4$.

  We proceed  to study the matrix equations  $X^r =\lambda SD^n,
  X^r =\lambda RD^n$ for an odd prime
  $r$, rational $\lambda$  and a natural $n$.
  Since for $r\neq 3$ the elements
  $S$ and $R$ have trivial $r$-th roots
 in $\mathcal L(t)$ it is enough to consider the case of $r=3$.
  By the argument in the proof of
  Proposition  \ref{Proposition20} we may further reduce
  our attention to $\lambda =1$ and $n=1$.
  By the assumption of circular primitivity
  the equations
  $X^3 =\lambda SD,  X^3 =\lambda RD$ have no solutions
  in $\mathcal L(t)$.

  It remains to observe that  the elements $Y$ and $Z$
  of order $6$ do not have square roots in $\mathcal G(t)$
  because $\det YD^n=\det ZD^n =12f^2$ and Proposition
  \ref{Proposition14} applies.

  The list of finite order elements in the Laxton group
  $\mathcal G(t)$ is completed by the multiplication of
  two elements of order $3$ and  three elements of order   $2$ to get
  four additional elements of order $6$.

  We have $(WS)^2 = 2f(2+t)D^{-1}R, (WR)^2 = -2f(2+t)D^{-1}S,
  (WY)^2 = 6f(2+t)D^{-1}S, (WZ)^2 = -6f(2+t)D^{-1}R$.
  \end{proof}

For a given $t = C_m(u)$, for some
natural $m$, and a rational primitive  $u$, we have the homomorphism
$\Xi_m:\mathcal G(t) \to \mathcal G(u)$. The kernel of this homomorphism
is a cyclic subgroup of order $m$, consisting of
the equivalence classes of the matrices
$\Phi_m(D_u^k), k =0,1,2,\dots, m-1$. The torsion subgroup of $\mathcal G(t)$
is the preimage of the torsion subgroup of $\mathcal G(u)$ under the
homomorphism $\Xi_m$, hence it is $m$ times larger.
Since $W_u^2 = D_u$ in $\mathcal L(u)$ then $\Phi_m(W_u)$ is
an element of order $2m$ in $\mathcal G(t)$. It follows that
the torsion subgroup of $\mathcal G(t)$ in the non-cyclotomic case
is isomorphic to $\mathbb Z_{2m} \times \mathbb Z_{2}$. For the same reason
in the cyclotomic  case
the torsion subgroup of $\mathcal G(t)$ is isomorphic
to $\mathbb Z_{2m} \times \mathbb Z_{4}$ if  $u$ is  circular and
circular primitive, and to
$\mathbb Z_{2m} \times \mathbb Z_{6}$
if $u$ is  cubic and cubic  primitive.

In terms of torsion sequences $\{x_n\}_{n\in\mathbb Z}$ of
$\mathcal G(u)$
the torsion elements in $\mathcal G(t)$ are the
subsequences $\{x_{km}\}_{k\in\mathbb Z}$.

Finally we consider a cubic and cubic primitive $u$, and
$t=C_{3m}(u)$ for  a natural $m$,
and  an associate value
$a$  of $t$, which by necessity is cubic but not cubic primitive.
%we have that $a = C_m(u')$ where $u'$ is the associate  value
%of  $C_{3^k}(u)$. It is clear that $u'$ is  primitive but not cubic primitive.
Due to our definitions $t':=C_3(t) = C_3(a)$.
  The sequence groups for all of these parameters are naturally isomorphic,
  and we choose to identify them with  $\mathcal L(t')$.
  Under such an identification we have $D_{t'}= D_t^{3}$ and
  $D_a=D_tS$.  The Laxton groups $\mathcal G(t)$ and $\mathcal G(a)$
  are factor groups $\mathcal G(t')/\{D_t\}$ and
  $\mathcal G(t')/\{D_a\}$, respectively. The torsion subgroup
  of $\mathcal G(t')$ is isomorphic to $\mathbb Z_{6m} \times \mathbb Z_{6}$.
  To get the torsion subgroups of   $\mathcal G(t)$ and $\mathcal G(a)$
  we factor $\mathbb Z_{6m} \times \mathbb Z_{6}$ by the elements
  $(3,0)$ and $(3,2)$, respectively. Hence we obtain that
the torsion subgroup of $\mathcal G(a)$ is isomorphic to
$\mathbb Z_{6m} \times \mathbb Z_{2}$, while it is isomorphic to
$\mathbb Z_{2m} \times \mathbb Z_{6}$ for
$\mathcal G(t)$.

In
the circular case of  a primitive and circular  primitive $u$,
and    $t=C_{2m}(u)$ and the associate value
$a$  of $t$, we obtain by a similar argument that
the   torsion subgroup of $\mathcal G(a)$ is isomorphic
to  $\mathbb Z_{4m} \times \mathbb Z_{2}$,
while the torsion subgroup of $\mathcal G(t)$ is isomorphic
to  $\mathbb Z_{2m} \times \mathbb Z_{4}$.

\section{The sets of prime divisors of recursive sequences}

Let us recall that for a fixed rational $t\neq 0,\pm 1,\pm 2$,
we denote by $\Pi_t$
the set of all odd primes with the exclusion of the divisors of the
denominators and numerators of $t$ and $\delta(t) = t^2-4$.
Not all of these exclusions are necessary for some of the forthcoming theorems.
However they are only finite, and they simplify the formulations.

%\in \mathcal L(t)
%A recursive sequence
%$[x_0 \ x_1]=X$ is called {\it reduced} if $x_0$ and $x_1$
%are integer and relatively prime. Clearly any element of $ \mathcal L(t)$
%has a reduced representative.

For a  sequence $X \in \mathcal L(t)$ with the reduced  representative
$[x_0 \ x_1]=X$,
we define the set $ \Gamma_X$ of prime divisors of $X$ by
\[
\Gamma_X= \{p\in \Pi_t \ | \ \exists\ n\in \mathbb Z, \ x_n =0  \mod  p\}.
\]

%of primes from $\Pi_t$ which divide
%some element of the sequence $\{x_n\}$, assuming that we choose
%the reduced representative to obtain $\{x_n\}$.

We have that   $p\in \Gamma_X$ if and only if there is an integer
$k$ such that $XD^k =x_{k+1} I  \mod  p$.
Let us recall that
if $X$ represents the sequence $\{x_k\}_{k\in \mathbb Z}$, then $X^{-1}$
represents  the sequence $\{x_{-k}\}_{k\in \mathbb Z}$.
It follows immediately that
 $\Gamma_X = \Gamma_{X^{-1}}$.  Similarly $\Gamma_X = \Gamma_{XD^{k}}$
for any integer $k$.

 We turn to some  general facts about the sets of prime divisors
of the recursive sequences.

\begin{proposition}
 \label{Proposition23}
  For any $X,Y \in \mathcal L(t)$, \
  $  \Gamma_X \cap \Gamma_Y \subset \Gamma_{XY}$.
  Further
\[
\Gamma_X \cap  \Gamma_{XY} = \Gamma_Y \cap \Gamma_{XY}=
\Gamma_X \cap \Gamma_Y.
  \]
\end{proposition}
\begin{proof}
  For $p\in\Gamma_X \cap \Gamma_{Y}$ there are integers
  $k$ and $l$ such that $XD^k =x_{k+1} I  \mod  p$ and
  $YD^l =y_{l+1} I  \mod  p$. Hence
   $XYD^{k+l} =x_{k+1}y_{l+1} I  \mod  p$, which proves the
  first part.

  To prove the second part we note that
  $\Gamma_X \cap  \Gamma_{XY} =\Gamma_{X^{-1}} \cap  \Gamma_{XY}
  \subset   \Gamma_{Y}$. The rest follows by the symmetric role
  played by $X$ and $Y$.
\end{proof}

%\begin{proposition}
%  For $X,Y \in \mathcal G(t)$,  $X= [x_0 \ x_1], Y= [y_0 \ y_1]$,
%  \[  \Gamma_X \cap \Gamma_Y \subset \Gamma_{XY}, \ \
%  \Gamma'_X \cap \Gamma_Y \subset \Gamma'_{XY},
%  \]
%  where $'$ denotes the complement of the respective set of primes
%  in $\Pi$.
%\end{proposition}
%\begin{proof}
%  We have $p\in \Gamma_X$ if and only if there is an integer
%  $k$ such that $XD^k =x_{k+1} I \ \mod \ p$. Hence
%   if  $p\in\Gamma_X \cap \Gamma_{Y}$ then there is also an
%  integer $l$ such that $YD^l =y_{l+1} I \ \mod \ p$.
%  It follows that $XYD^{k+l} =x_{k+1}y_{l+1} I$, which proves the
%  first part.
%
%  To prove the second part we assume that $p\notin\Gamma_{X}$ and
%   $p\in\Gamma_Y \cap \Gamma_{XY}$.
%  Hence $p\in\Gamma_{Y^{-1}} \cap \Gamma_{XY}$ and by the first part
%  $p\in \Gamma_{X}$, a contradiction which ends the proof.
%  \end{proof}
%\begin{corollary}
%  For any  $X,Z\in \mathcal G(t) $ such that  $ Z^2 = \lambda I$,
%  for some scalar $\lambda$,
%  \[  \Gamma_{X}\cap  \Gamma_{Z} = \Gamma_{X}\cap  \Gamma_{ZX} = \Gamma_{ZX}\ca%p  \Gamma_{Z}.
%  \]
%\end{corollary}
%\begin{proof}
%  We get immediately
%   $\Gamma_{X}\cap  \Gamma_{Z} \subset  \Gamma_{ZX}$ and
%  $\Gamma_{ZX}\cap  \Gamma_{Z} \subset  \Gamma_{X}$.
%Further
%$\Gamma_{X}\cap  \Gamma_{ZX}=
%\Gamma_{X^{-1}}\cap  \Gamma_{ZX}\subset  \Gamma_{Z}$.
%\end{proof}
%Clearly we can apply this Corollary to the elements
%$C,W,V$.

It follows directly from Proposition \ref{Proposition61} that $\Gamma_{X^2}$
is the union of the following four subsets.
\begin{proposition}
 \label{Proposition24}
  For any  $X\in \mathcal L(t)$,
  \[  \Gamma_{X^2} =   \Gamma_{X} \cup  \Gamma_{CX} \cup
  \Gamma_{WX} \cup  \Gamma_{VX}.
  \]
\end{proposition}
We give another proof which makes no reference to the factorizations
of Proposition \ref{Proposition61}.
\begin{proof}
  Let us note that for any $p\in \Pi_t$ the element  $X^2$ in
  the cyclic group $\mathcal L_p(t)$
  has exactly two square roots, $X$ and $CX$. Further the element
  $D^{-1}$ also has exactly two square roots, $W$ and $V$.

     $p\in \Gamma_{X^2}$ if and only if  there is an integer $k$ such that
  $X^2 =  D^k$ in $\mathcal L_p(t)$.

  If $k=2l$  then $\left(D^{l}\right)^2 = X^2$, and hence in the
  cyclic group $\mathcal L_p(t)$ the element
  $D^l$ is a square root of $X^2$. It follows immediately that
  either $D^l = X$ or $D^l= CX$ in $\mathcal L_p(t)$.

  If $k=2l-1$  then $\left(XD^{-l}\right)^2 =
   D^{-1}$ in $\mathcal L_p(t)$
  and, since we know two square roots of $D^{-1}$ in $\mathcal L(t)$,
  we conclude that
  $XD^{-l}= W$ or
  $XD^{-l} = V$ in $\mathcal L_p(t)$.
\end{proof}

\begin{figure}[htbp]
\includegraphics[width=40mm]{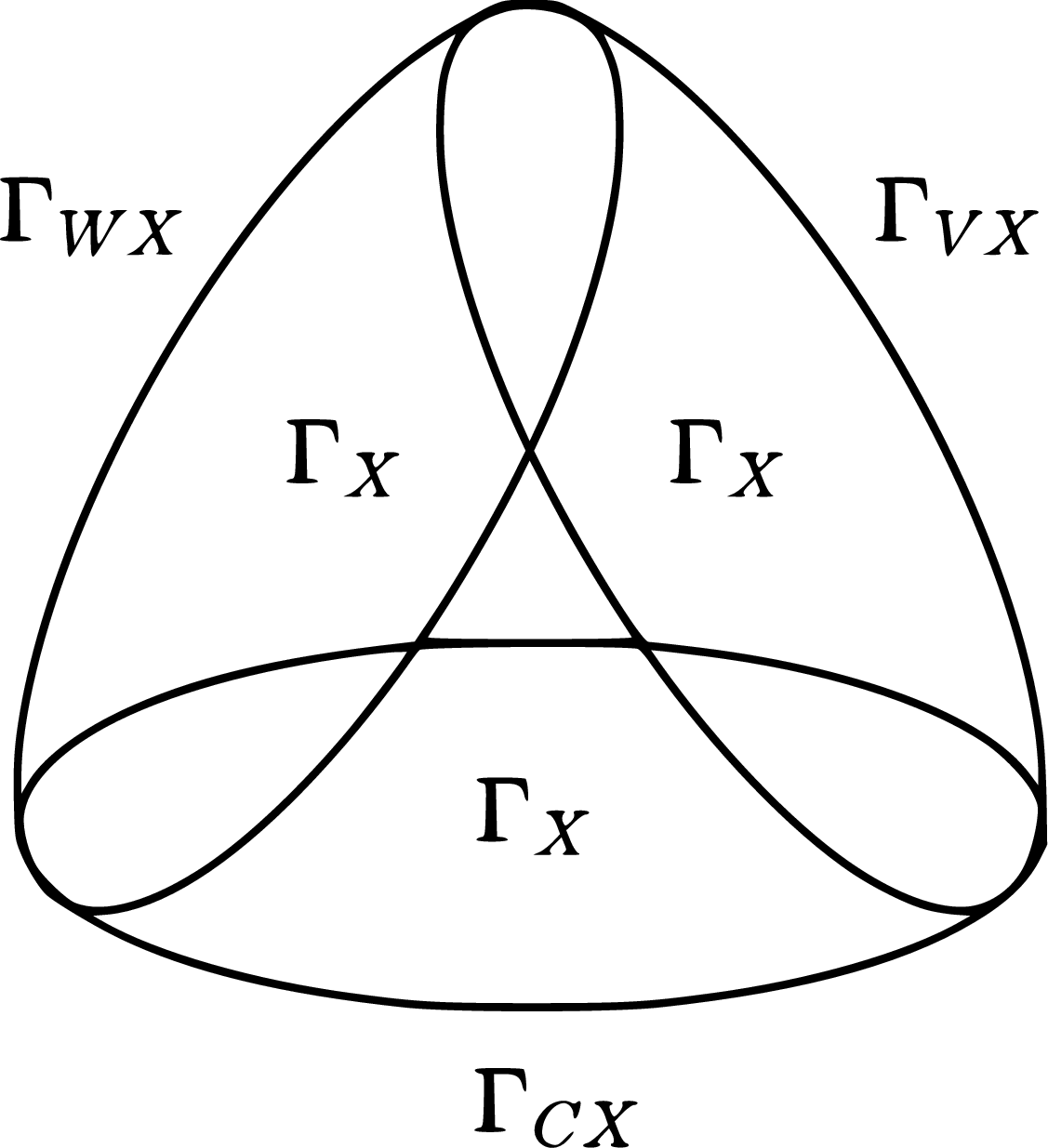}
\caption{Sets of prime divisors of $X, WX, VX, CX$.}
\label{fig1}
\end{figure}

The four sets in the last Proposition are not disjoint.
Actually the sum of any three of them will give the whole
$\Gamma_{X^2}$.
In general four sets produce a partition into $16$ subsets.
It turns out that in our case only $6$ elements of
the partition may be nonempty, Figure \ref{fig1}.

%=======================================================
%=======================================================

\begin{proposition}
 \label{Proposition25}
  For any  $X\in \mathcal L(t)$ the following six
  sets are disjoint
  \[
      \Gamma_{X} \cap  \Gamma_{CX}, \
\Gamma_{X} \cap  \Gamma_{WX},  \
  \Gamma_{X} \cap  \Gamma_{VX}, \
  \Gamma_{WX} \cap  \Gamma_{VX}, \
\Gamma_{CX} \cap  \Gamma_{WX},  \
  \Gamma_{CX} \cap  \Gamma_{VX},
  \]
  and their union is all of   $\Gamma_{X^2}$.
  Further
  \[
  \begin{aligned}
   &   \Gamma_{X} \cap  \Gamma_{CX}=\Gamma_{C} \cap  \Gamma_{CX}, \
\Gamma_{X} \cap  \Gamma_{WX}=\Gamma_{W} \cap  \Gamma_{WX},  \
  \Gamma_{X} \cap  \Gamma_{VX}=  \Gamma_{V} \cap  \Gamma_{VX},\\
  &  \Gamma_{WX} \cap  \Gamma_{VX} \subset \Gamma_{C},\ \ \
\Gamma_{CX} \cap  \Gamma_{WX}\subset \Gamma_{V},  \ \ \
\Gamma_{CX} \cap  \Gamma_{VX}\subset \Gamma_{W}.
\end{aligned}
  \]
\end{proposition}
In the proof we need the following Proposition.
\begin{proposition}
 \label{Proposition26}
  $\Gamma_{W},\Gamma_{V}$ and   $\Gamma_{C}$ are disjoint
  and their union is all of $\Pi_t$.
\end{proposition}
It was established in \cite{Wo}, we will give an independent proof
later on.
\begin{proof}
To check that the union is all of $\Gamma_{X^2}$ we observe that by
  Proposition \ref{Proposition23} we have the following equalities
  \[
  \begin{aligned}
   &   \Gamma_{X} \cap  \Gamma_{CX}=\Gamma_{C} \cap  \Gamma_{X}, \
\Gamma_{X} \cap  \Gamma_{WX}=\Gamma_{W} \cap  \Gamma_{X},  \
\Gamma_{X} \cap  \Gamma_{VX}=  \Gamma_{V} \cap  \Gamma_{X}.
  \end{aligned}
  \]
  By Proposition \ref{Proposition26} we conclude that
   these sets partition $\Gamma_X$.
  Replacing $X$ successively with $CX,WX $
  and $VX$ we will arrive at their partitions as well.
  One needs to recall that  $WV = CD^{-1}$,
   $WC = V, VC =W$ in $\mathcal L(t)$.
For example replacing $X$ with $VX$ we get
\[
  \begin{aligned}
    &   \Gamma_{VX} \cap  \Gamma_{CVX}=\Gamma_{VX} \cap  \Gamma_{WX}=
    \Gamma_{C} \cap  \Gamma_{VX}, \\
&\Gamma_{VX} \cap  \Gamma_{WVX}=\Gamma_{VX} \cap  \Gamma_{CX}=\Gamma_{W} \cap  \Gamma_{VX},  \\
    &\Gamma_{VX} \cap  \Gamma_{V^2X}= \Gamma_{VX} \cap  \Gamma_{X}=
    \Gamma_{V} \cap  \Gamma_{VX}.
  \end{aligned}
  \]
\end{proof}

Proposition \ref{Proposition26} generalizes to the following Corollary.
\begin{corollary}
 \label{Corollary27}
  If $\Gamma_{X^2} = \Gamma_{X}$ then the sets
  $\Gamma_{CX},\Gamma_{WX}$ and   $\Gamma_{VX}$ are disjoint
  and their union is all of $\Gamma_{X^2}$. Moreover
  \[
 \Gamma_{CX}=\Gamma_{C} \cap  \Gamma_{X}, \
\Gamma_{WX}=\Gamma_{W} \cap  \Gamma_{X},  \
 \Gamma_{VX}=  \Gamma_{V} \cap  \Gamma_{X}.
\]
\end{corollary}
\begin{proof}
  It follows directly from  Proposition \ref{Proposition25} when we observe that
  \[
  \Gamma_{X^2} \setminus \Gamma_{X} =
  \left(\Gamma_{WX} \cap  \Gamma_{VX}\right)\cup
\left(\Gamma_{CX} \cap  \Gamma_{WX}\right)  \cup
\left(\Gamma_{CX} \cap  \Gamma_{VX}\right).
\]
\end{proof}
%Choosing  $X=I$ we obtain the following Corollary.
In the cubic case $S^2=S^{-1}=R$, and $CS=Y, VS = WZ$ in $\mathcal L(t)$
which gives us yet another Corollary.
\begin{corollary}
 \label{Corollary28}
  In the  cubic case   $\Gamma_{WS},\Gamma_{Y}$ and   $\Gamma_{WY}$ are disjoint
  and their union is all of $\Gamma_{S}$.
  \end{corollary}
Ballot, \cite{B}, established, under some genericity conditions,
that the prime density of $\Gamma_{S}$ is $3/4$ and for
$\Gamma_{Y} \cup\Gamma_{WY}$ it is $1/2$. Assuming these genericity
conditions for both $t$ and the twin parameter $-t$ it follows,
using the last Corollary and the recombination of Propositions
\ref{Proposition11} and \ref{Proposition12},
that the prime densities of the three disjoint sets $\Gamma_{WS},\Gamma_{Y}$
and   $\Gamma_{WY}$ are each equal to $1/4$.
Indeed the sequences $WS$ in $\mathcal L(t)$ and
$VS = WY$ in $\mathcal L(-t)$ have essentially the same primes divisors,
and hence
$\Gamma_{WS}$ and   $\Gamma_{WY}$ have the same prime densities.
The same claim can be made about
$VS = WY$ in $\mathcal L(t)$
and $WS$ in $\mathcal L(-t)$.
It follows by simple arithmetic that the three sets have equal
prime densities.

 The sets of divisors can be described in terms of the groups
$\mathcal L_p(t)$, which are cyclic by   Theorem \ref{Theorem16}.
 For an element $X$ of the sequence group $\mathcal L(t)$
 let us denote by $ord_p(X)$ the order
 of $X$
in the group $\mathcal L_p(t)$, if $\det \ X \neq 0  \mod  p$.

 The following Proposition
is essentially equivalent to the theorem  of Hall, \cite{H2}.
Let us observe   first that if for a reduced sequence $\det \ X = 0 \mod  p$
then  $p\notin \Gamma_X$. Indeed, since $\det \ X=x_1^2-tx_1x_0+x_0^2$,
in such a case if $p \ | \ x_0$ then
$p \ | \ x_1$, contradicting the assumption that the sequence is reduced.
Further, by the same argument, for the shifted sequence
$XD^k$, for any integer $k$, if $x_k = 0  \mod  p$ then
$x_{k+1} = 0  \mod  p$. But then all elements of the sequence
vanish $ \mod  p$, which is again the contradiction.
\begin{proposition}{}
 \label{Proposition29}
  For any $X\in \mathcal L(t)$
  \[
  \Gamma_X =\{p\in \Pi_t\ | \  ord_p(X) \ | \ ord_p(D_t)\}.
  \]
  \end{proposition}
%\begin{proposition}
%  For any $p\in \Pi(t)$, we have
%  $p\in\Gamma_X$ if and only if $ord(X) \ | \ ord(D_t)$,
%  where $X$ and $D_t$ are considered as the elements of the cyclic
%  group $\mathcal L_p(t)$.
%  \end{proposition}
\begin{proof}
  As noted before $p\in \Gamma_X$ if and only if there is
  an integer $k$ such that $X=x_{k+1}D^{-k}  \mod  p$.
  The last equality can be considered as the equality in the
  cyclic group $\mathcal L_p(t)$.
  Now our Proposition
  is a consequence of the following property of cyclic groups.
  For two  elements  $x,y$ of a finite cyclic group
  $ord(x) \ | \ ord(y)$ if and only if
 $x$ is contained in the subgroup  generated
  by $y$.
\end{proof}

As an application of this criterion we prove the following
more detailed version of Proposition \ref{Proposition26}
and Corollary \ref{Corollary28}. Let us define the
{\it index of appearance  $\xi= \xi(t,p)$} $=ord_p(W)$.
This definition is justified by the fact that $\xi(t,p)$ is
actually the classical  index of appearance for the
 Lucas sequence $\{L_n\}_{n\in\mathbb Z}$
 for the two parameters $(T,Q)$, i.e., it is the smallest natural index
 $n$ such that $L_n=0$. Indeed, let us  denote $W =\{w_{2n-1}\}_{n\in\mathbb Z}$
 and $D =\{u_{n+1}\}_{n\in\mathbb Z}$. If $\xi = 2k-1$ then
 $W^\xi =I = W^{-\xi} = D^{k}W$ which leads by formula (2) to
 $0=w_\xi =Q^{1-k} L_\xi$.  If $\xi = 2k$ then
 $W^\xi =I = W^{-\xi} = D^{k}$ which leads by formula (2) to
 $0=u_k =T^{-1}Q^{1-k} L_\xi$. The rest of the argument is routine.
\begin{proposition}
  \label{Proposition30}
  For any $t\neq 0,\pm 1, \pm 2$
 \[
  \begin{aligned}
  \Gamma_{W} = \{p\in\Pi_t\ |\ \xi = 1\mod 2\},
 \ & \Gamma_{V} = \{p\in\Pi_t\ |\ \xi = 2\mod 4\},\\
 & \Gamma_{C} = \{p\in\Pi_t\ |\ \xi = 0\mod 4\}.
\end{aligned}
  \]
  For any cubic $t$
 \[
  \begin{aligned}
  \Gamma_{WS} = \{p\in\Pi_t\ |\ \xi = 3\mod 6\},
 \ & \Gamma_{VS} = \{p\in\Pi_t\ |\ \xi = 6\mod 12\},\\
  & \Gamma_{CS} = \{p\in\Pi_t\ |\ \xi = 0\mod 12\}.
\end{aligned}
 \]
   \end{proposition}
The first part  appeared as ``trichotomy'' in \cite{B-E}.
It was also proven in \cite{Wo}, where  the sets
are denoted as $ \Pi_0, \Pi_1, \Pi_*$.
The second part was essentially proven in \cite{B}.
   \begin{proof}
     The cyclic group $\mathcal L_p(t)$ is isomorphic to the additive
     group $\mathbb Z_N$
     for the appropriate $N=p\pm 1$. We do not have any canonical
     isomorphism of the two,
     nevertheless we find it convenient to do calculations
     in $\mathbb Z_N$.

     Since $C$ is of order $2$ in $\mathcal L_p$ then  $C = \frac{N}{2}
     \in \mathbb Z_N$. Let $W =w \in \mathbb Z_N$.
     Since we do calculations in a cyclic group without loss of
     generality we can assume that
     $\xi = ord_p(W) = \frac{N}{w}$.
     Since $W^2 = D^{-1} = V^2$
     in $\mathcal L_p$ then $D^{-1} = 2w\in \mathbb Z_N$
     and $V = w+\frac{N}{2}\in \mathbb Z_N$.

We get then that $ord_p(D^{-1}) = ord_p(D)$ is equal to $\xi$, if $\xi$ is odd,
and to $\xi/2$, if $\xi$ is even. It follows that $p \in \Gamma_W$ if and
only if $\xi$ is odd.

Further $p\in \Gamma_C$ if and only if $2 = ord_p(C)$
divides  $ord_p(D)$, which is     $\xi$ or $\xi/2$.
It follows that
$p\in \Gamma_C$ if and only if
     $\xi = 0  \mod  4$.

Since $V= \frac{1}{2}w\left(2+\xi\right)$ we get that if $\xi$ is odd then
$ord_p(V)=2\xi$ and $p\notin\Gamma_V$.
If
$\xi$ is even and $1+\frac{\xi}{2}$ is odd then
$ord_p(V)=\xi$ and $p\notin\Gamma_V$.
If
$\xi$ is even
and $1+\frac{\xi}{2}$ is even than      $ord_p(V)=\xi/2$
and $p\in\Gamma_V$. We can now conclude that $p\in\Gamma_V$
if and only if $\xi = 2  \mod  4$.

In the cubic case $S= \frac{N}{3}$ and $WS = \frac{1}{3}w(3+\xi)$.
If $3\nmid \xi$ then $ord_p(WS)=3\xi$ and  hence $p\notin \Gamma_{WS}$.
If $3\nmid \xi$ then $ord_p(WS)=\xi$ and hence
$p\in \Gamma_{WS}$ if and only if $\xi$ is odd.

We have  $CS = \frac{N}{2}+\frac{N}{3}$ and hence
$ord_p(CS)= 6$. Since $ord_p(D)$ is $\xi$ or $\xi/2$ then
clearly $p\in \Gamma_{CS}$ if and only if
     $\xi = 0  \mod  12$.

Finally $VS= w+\frac{N}{2}+\frac{N}{3} =\frac{1}{6}w\left(6+5\xi\right)$.
It follows that if $\xi$ is odd then $ord_p(VS)$ must be even
and hence it cannot divide $ord_p(D)=\xi$. If $3\nmid \xi$ then
$3\mid ord_p(VS)$ and then also $p\notin \Gamma_{VS}$. If
$6\mid \xi$ then   $VS= w\left(1+5\frac{\xi}{6}\right)$ and
$ord_p(VS) =\xi$ if $\xi/6$ is even, and $ord_p(VS) =\xi/2$
if $\xi/6$ is odd. Hence
$p\in \Gamma_{VS}$ if and only if
     $\xi = 6  \mod  12$.
   \end{proof}

% ----------------------------------------------------->
\section{Density estimates and Independence Conjecture}
% ----------------------------------------------------->
\begin{theorem}
 \label{Theorem31}
 For any sequence $X\in \mathcal R(t)$ the set of divisors $\Gamma_X$
 is contained in the subset of primes in $\Pi_t$ such that
 $\det X$ is a square $\mod p$.
  
If $\det X$ is not a rational square then
the upper prime density of $\Gamma_X $ does not exceed $1/2$.
\end{theorem}
%---
\begin{proof}
If $p \in \Gamma_{X}$ then there is integer $k$ such that $X D^k = x_{k+1} I {\rm \;\mod \;} p$.
Taking determinants of both sides
we obtain $\det X = x^2_{k+1} {\rm \;\mod \;} p$.
\end{proof}
%--
It follows from this Theorem that for two sequences
$X, Y \in \mathcal{L}(t)$, we have the following table,
where we indicate the location of the elements of $\Gamma_X$ and $\Gamma_Y$.
\begin{table}[htbp]
  \caption{The sets of divisors for two sequences $X, Y \in \mathcal{L}(t)$}
$$\begin{array}{|c|c|c|}
\hline
                    & \det X = \square        & \det X \neq \square \\
 \hline
\det Y = \square    & \Gamma_X \cup \Gamma_Y &  \Gamma_Y\\
\hline
\det Y \neq \square & \Gamma_X                  &  \emptyset\\
 \hline
\end{array}$$
\end{table}
In particular by the Frobenius Density Theorem, \cite{S-L},
if $\det X, \det Y$ and $\det XY$ are not rational
squares then the upper prime densities of
$\Gamma_X \cap \Gamma_Y$ and $\Gamma_X \cup \Gamma_Y$ do not
exceed $1/4$ and $3/4$, respectively.

For a sequence $[\tilde{x}_0 \;\tilde{x}_1] = \widetilde{X} \in \mathcal{L}(T, Q)$, with simple integer parameters $T, Q$, we consider the sets of prime divisors
of the even numbered elements, and the odd numbered elements. As observed previously for the sequence
$[Q \tilde{x}_0\; \tilde{x}_2]  = X \in \mathcal{L}(t)$, $t =\frac{T^2}{Q} - 2$, the sets $\Gamma_X$ and $\Gamma_{WX}$ are the sets of
divisors of even and odd numbered elements of $\widetilde{X}$, respectively.

We have $\det X = T^2  \det \widetilde{X}$ and
$\det W = 2 + t = \frac{T^2}{Q} $.
Under the assumption that $\det \widetilde{X}$, $Q$ and $Q \det \widetilde{X}$ are not rational squares
we get an approximate “independence” of the sets of divisors, which can be
illustrated in the following table.

\begin{table}[htbp]
  \caption{The sets of divisors for the sequences $X, WX$}
$$\begin{array}{|c|c|c|}
\hline
                & \det \widetilde{X} = \square            & \det \widetilde{X} \neq \square \\
 \hline
Q = \square     & \Gamma_X\cup \Gamma_{WX} &  \emptyset\\
\hline
Q \neq \square  & \Gamma_{X}                  &  \Gamma_{WX}\\
 \hline
\end{array}$$
\end{table}

More precisely, as observed above the upper prime densities of $\Gamma_{X}$, $\Gamma_{WX}$ and
$\Gamma_{X} \cap \Gamma_{WX}$ do not exceed $1/2$ and $1/4$, respectively.

Schinzel, \cite{Sch}, considered the sequence $B =\{a^n-b\}_{n\in \mathbb N}$
for integer
$b,a\neq 1$, and proved that if $b$ is not a power of $a$
then the set of primes that do not divide any element of the
sequence is infinite. The sequence
$B = [1-b\; a-b] \in \mathcal R(a+1, a)$
so that we can apply our estimate on the size of the sets of prime divisors.
We get that if $a, \ \det B = b(a -1)^2$ and $ab$ are not rational
squares then the set of primes not dividing any element of  $B$
contains a subset of prime density $1/4$. This set is equal to the set
of primes $p$ such that $b$ is a square $\mod p$ and $a$ is not a square
$\mod p$. Let us note that  such a result has a simple direct proof based
on the analysis of the equations $a^k = b \mod p$ separately for
even and odd values of $k$. Finally let us observe that under the
same conditions all the sets
of prime non-divisors for odd powers $(a-1)^{1-k}B^k= \{a^n-b^k\}_{n\in \mathbb N}$
contain the common subset of density $1/4$.  It is clear that
the set of prime non-divisors for the power $B^k$ is not smaller
than for $B^{kl}, l > 1$ (it follows also from
Proposition \ref{Proposition29}).

While we cannot apply  Theorem \ref{Theorem31}  directly to sequences
$X^2$ for $X \in \mathcal R(t)$, we get a weaker estimate using
Proposition \ref{Proposition24}.
\begin{proposition}
 \label{Proposition32}
 For any sequence $X\in \mathcal R(t)$
 the set of non-divisors $\Pi_t\setminus \Gamma_{X^2}$
 contains  the subset of primes $p$ in $\Pi_t$ such that
  $\det X$ is not a square $\mod p$ and
 $\det W$ and $\det V$ are squares $\mod p$.

 If for $a=\det X,\ b=\det W=2+t,\ c=\det V= 2-t$ the degree of the field
 extension $[\mathbb Q(\sqrt{a},\sqrt{b},\sqrt{c}) : \mathbb Q] =8$
 then the upper prime density of $\Gamma_{X^2}$  does not exceed $7/8$.
\end{proposition}
\begin{proof}
  By Proposition 25 we have $\Gamma_{X^2} = \Gamma_{X}\cup
  \Gamma_{WX}\cup \Gamma_{VX}$. We apply now Theorem \ref{Theorem31}
  to the three sequences $X, WX$ and $VX$.

  The second part follows from the Frobenius Density Theorem.
\end{proof}
Note that under the assumptions of this Proposition the sets of non-divisors
of $X^2$ and $WX^2$ contain the  subsets of primes of density $1/8$
and $1/2$ respectively, but they are disjoint (in the former
$\det W$ is a square $\mod p$ and in the latter it is not a square
$\mod p$).

For a subset of odd primes $A$ we denote by $|A|$ the prime density
of the set, if it exists. We say that a parameter $t$ is {\it generic}
if it is twin-primitive and $r$-primitive for any prime $r\geq 3$,
and it is not cyclotomic (i.e., circular or cubic). This condition can be
effectively checked. Indeed if $x = \frac{a}{b}$
with relatively prime integers $a, b \geq 1$, then the reduced
fraction $C_n(x)$ has $b^n$ in the denominator. Hence the necessary
condition for the equation $C_r(x)=t$ to have a rational solution
is that the denominator
of $t$ is an $r$-th power (in the reduced form).
Once we know which $r$ to consider we can find all the rational
solutions. It remains to address the case of integer $t, |t|\geq 3$.
It is clear then that if there is a rational solution $x_0$ of $C_r(x)=t$
then it must be integer and $|x_0|\geq 3$.
The polynomial $C_r(x)$ is increasing for $|x|> 2$. Hence
the only possible values of $r$ are such that $C_r(3) \leq |C_r(x_0)| =|t|$.

{ \bf Conjecture on the exact independence of even and odd
  numbered elements of most
  recursive sequences}

{ \it For any generic $t$ and a sequence $X$  of infinite order
  in the Laxton group $\mathcal{G}(t)$ which is not any power in that group,
  if $\det X, \det (WX),\det (CX) $ and $\det (VX)$ are not rational squares
  then the sets of
  primes $\Gamma_{X}$ and $\Gamma_{WX}$ have prime densities and
are ``independent''}
$$|\Gamma_{X} \cap \Gamma_{WX}| = |\Gamma_{X}| |\Gamma_{WX}|.$$
%
% \newpage

While we believe that most sequences are not any power in the Laxton group
$\mathcal G(t)$,
we do not know an  effective way to check this condition in examples.

If the conjecture is valid then it applies equally well
to $CX$ and hence we get also the independence of
$\Gamma_{CX}$ and $\Gamma_{VX}$. Further we also get the independence
for the respective sequences in $\mathcal L(-t)$. Hence by the recombination
of Proposition \ref{Proposition7} we would obtain that independence applies
also to the pairs of sets
$\Gamma_{X}, \Gamma_{VX}$ and
$\Gamma_{CX},\Gamma_{WX}$. It turns out that the independence of the four
pair of sets implies that
\[
|\Gamma_{X} \cap \Gamma_{W}| = |\Gamma_{X} \cap \Gamma_{V}| =
|\Gamma_{CX} \cap \Gamma_{WX}| = |\Gamma_{CX} \cap \Gamma_{VX}|,
\]
and in view of Table 1 this leads to
\[
|\Gamma_{X}|=|\Gamma_{CX}|, \ |\Gamma_{WX}|=|\Gamma_{VX}|.
\]
To prove it let us denote the prime densities of the six subsets
shown in Table 1 as
\[
\begin{aligned}
  a &=  |\Gamma_{X} \cap \Gamma_{W}|, b = |\Gamma_{X} \cap \Gamma_{V}|,
  c = |\Gamma_{X} \cap \Gamma_{CX}|,\\
\alpha &= |\Gamma_{CX} \cap \Gamma_{VX}|, \beta = |\Gamma_{CX} \cap \Gamma_{WX}|,
\gamma = |\Gamma_{WX} \cap \Gamma_{VX}|.
\end{aligned}
\]
The independence of the four pairs
$(\Gamma_{X}, \Gamma_{WX})$,
$(\Gamma_{X},\Gamma_{VX})$, $(\Gamma_{CX}, \Gamma_{VX})$ and
$(\Gamma_{CX},\Gamma_{WX})$ give us the equalities
\[
\begin{aligned}
&(1)\ \ \ (a+b+c)(a+\beta+\gamma)  = a,\\
&(2)\ \ \   (a+b+c)(\alpha+b+\gamma)  = b,\\
&(3) \ \ \ (\alpha+\beta+c)(\alpha+b+\gamma)  = \alpha,\\
&(4) \ \ \  (\alpha+\beta+c)(a+\beta+\gamma)  = \beta.
\end{aligned}
  \]
Dividing (1) by (2) we obtain
$
 b(a+\beta+\gamma) = a(\alpha+b+\gamma),
$
which yields
\[
(5) \ \ \  b(\beta+\gamma) = a(\alpha+\gamma).
\]
Dividing (3) by (4) we get
$
 \beta(\alpha+b+\gamma) = \alpha(a+\beta+ \gamma),
$
which leads to
\[
(6)\ \ \ \beta(b+\gamma) = \alpha(a+ \gamma).
\]
Subtracting  (6) from (5)
gives us
$
\gamma(\beta-b) = \gamma(\alpha-a),
$
or equivalently
\[
(7) \ \ \   a-b = \alpha - \beta.
\]
Dividing (1) by (4)  we get
$
\beta(a+b+c) = a(\alpha+\beta+c)
$
or
\[
(8) \ \ \  \beta(b+c) = a(\alpha+c).
\]
Dividing the (2) by (3)
delivers
$
b(\alpha+\beta+c) = \alpha(a+b+ c)
$
or
\[
(9)   \ \ \ b(\beta+c) = \alpha(a+ c).
\]
Subtracting (8) from (9) produces
$
c(b-\beta) = c(\alpha-a)
$
or
\[
(10) \ \ \    a + b = \alpha  + \beta.
\]
Now
(7)  and   (10)  imply
\[
(11) \ a = \alpha, \ \ \ \ b = \beta.
\]
Substituting (11) into (1) and (2) we conclude that
$a = b$, so that finally
$a = b = \alpha = \beta$.

Numerical experiments support the conjecture, which is shown in the following
Table 3.
%Numerical experiments for the first $1200$ primes, Table \ref{tab1}.
%----

%=======================================================
\begin{table}[htbp]
\caption{The prime density of $|\Gamma_{X}|$ calculated for the first $1200$ primes}
\begin{center}
\begin{tabular}{cccccc}
\hline
%\textbf{Table}&\multicolumn{3}{|c|}{\textbf{Table Column Head}} \\
%\cline{2-4}
%\textbf{}(T,Q) & \textbf{\textit{$[x_0 \;x_1]$}} & \textbf{\textit{$|\Gamma_{X_{{\rm even}}}|$}}& \textbf{\textit{$|\Gamma_{X_{{\rm odd}}}|$}} & \textbf{\textit{$|\Gamma_{X}|$}}  & \textbf{\textit{$|\Gamma_{X_{{\rm even}}} \cap \Gamma_{X_{{\rm odd}}}|$}}\\
%\textbf{}(T,Q) & \textbf{\textit{$[x_0 \;x_1]$}} & \textbf{\textit{$|\Gamma_{\widetilde{X}_{{\rm even}}}|$}}& \textbf{\textit{$|\Gamma_{\widetilde{X}_{{\rm odd}}}|$}} & \textbf{\textit{$|\Gamma_{\widetilde{X}}|$}}  & \textbf{\textit{$|\Gamma_{\widetilde{X}_{{\rm even}}} \cap \Gamma_{\widetilde{X}_{{\rm odd}}}|$}}\\
%\hline
\textbf{}(T,Q) & \textbf{\textit{$[\widetilde x_0 \;\widetilde x_1]$}} & \textbf{\textit{$|\Gamma_{X}|$}}& \textbf{\textit{$|\Gamma_{WX}|$}} &  \textbf{\textit{$|\Gamma_{X} \cap \Gamma_{WX}|$}} & \textbf{\textit{$|\Gamma_{X}| \; |\Gamma_{WX}|$}} \\
\hline
        (5,3)  &      [17 $\;$11]  &          0.356           &      0.356       &           0.126      &   0.127\\
%\hline
        (4,11) &       [3 $\;$8]   &          0.328           &     0.334        &           0.111       &   0.110 \\
%\hline
         (3,7) &       [2 $\;$5]   &          0.357           &     0.355        &           0.123       &  0.127 \\
%\hline
       (3,-2)  &       [4 $\;$15]  &        0.353             &     0.340        &           0.116       &   0.120       \\
%\hline
        (7,11) &       [3 $\;$2]   &        0.340             &    0.340         &        0.115           & 0.116         \\
%\hline
       (2,-5)  &       [3 $\;$14]  &         0.343             &    0.339           &     0.111               &   0.116       \\
\hline
% \multicolumn{4}{l}{ $^{\mathrm{a}}$ Sample of a Table footnote.}
\end{tabular}
\label{tab1}
\end{center}
\end{table}
%=======================================================
We also inspected  the sets of divisors for squares (Table 4) and third powers
(Table 5) of $X$.
%=======================================================
\begin{table}[htbp]
\caption{The prime density of $|\Gamma_{X^2}|$ calculated for the first $1200$ primes}
\begin{center}
\begin{tabular}{cccccc}
\hline
%\textbf{Table}&\multicolumn{3}{|c|}{\textbf{Table Column Head}} \\
%\cline{2-4}
%\textbf{}(T,Q) & \textbf{\textit{$[x_0 \;x_1]$}} & \textbf{\textit{$|\Gamma_{X_{{\rm even}}}|$}}& \textbf{\textit{$|\Gamma_{X_{{\rm odd}}}|$}} & \textbf{\textit{$|\Gamma_{X}|$}}  & \textbf{\textit{$|\Gamma_{X_{{\rm even}}} \cap \Gamma_{X_{{\rm odd}}}|$}}\\
%\textbf{}(T,Q) & \textbf{\textit{$[x_0 \;x_1]$}} & \textbf{\textit{$|\Gamma_{\widetilde{X}_{{\rm even}}}|$}}& \textbf{\textit{$|\Gamma_{\widetilde{X}_{{\rm odd}}}|$}} & \textbf{\textit{$|\Gamma_{\widetilde{X}}|$}}  & \textbf{\textit{$|\Gamma_{\widetilde{X}_{{\rm even}}} \cap \Gamma_{\widetilde{X}_{{\rm odd}}}|$}}\\
\textbf{}(T,Q) &
& \textbf{\textit{$|\Gamma_{X^2}|$}}
& \textbf{\textit{$|\Gamma_{WX^2}|$}}
&  \textbf{\textit{$|\Gamma_{X^2} \cap \Gamma_{WX^2}|$}}
& \textbf{\textit{$|\Gamma_{X^2}| \; |\Gamma_{WX^2}|$}} \\
\hline
        (5,3)  &   [-1071 $\;$-746] &    0.722269           &  0.316931         &    0.261885      &  0.228909 \\
%\hline
        (4,11) &   [12 $\;$-35]     &    0.670559           &    0.289408       &    0.231026      & 0.194065 \\
%\hline
         (3,7) &    [8 $\;$-3]       &     0.70392          &     0.297748      &     0.237698     &  0.209591 \\
%\hline
       (3,-2)  &   [72 $\;$257]     &    0.692244           &    0.296914       &    0.230192      &  0.205537\\
%\hline
        (7,11) &   [-51 $\;$-95]    &   0.688073            &    0.30442       &    0.237698      &  0.209463\\
%\hline
       (2,-5)  &  [66 $\;$241]      &  0.686405            &   0.297748        &    0.226856      & 0.204376 \\
\hline
% \multicolumn{4}{l}{ $^{\mathrm{a}}$ Sample of a Table footnote.}
\end{tabular}
\label{tab1}
\end{center}
\end{table}
%=======================================================

%=======================================================
\begin{table}[htbp]
\caption{The prime density of $|\Gamma_{X^3}|$ calculated for the first $1200$ primes}
\begin{center}
\begin{tabular}{cccccc}
\hline
%\textbf{Table}&\multicolumn{3}{|c|}{\textbf{Table Column Head}} \\
%\cline{2-4}
%\textbf{}(T,Q) & \textbf{\textit{$[x_0 \;x_1]$}} & \textbf{\textit{$|\Gamma_{X_{{\rm even}}}|$}}& \textbf{\textit{$|\Gamma_{X_{{\rm odd}}}|$}} & \textbf{\textit{$|\Gamma_{X}|$}}  & \textbf{\textit{$|\Gamma_{X_{{\rm even}}} \cap \Gamma_{X_{{\rm odd}}}|$}}\\
%\textbf{}(T,Q) & \textbf{\textit{$[x_0 \;x_1]$}} & \textbf{\textit{$|\Gamma_{\widetilde{X}_{{\rm even}}}|$}}& \textbf{\textit{$|\Gamma_{\widetilde{X}_{{\rm odd}}}|$}} & \textbf{\textit{$|\Gamma_{\widetilde{X}}|$}}  & \textbf{\textit{$|\Gamma_{\widetilde{X}_{{\rm even}}} \cap \Gamma_{\widetilde{X}_{{\rm odd}}}|$}}\\
%\hline
\textbf{}(T,Q) &
& \textbf{\textit{$|\Gamma_{X^3}|$}}&
\textbf{\textit{$|\Gamma_{WX^3}|$}}
&  \textbf{\textit{$|\Gamma_{X^3} \cap \Gamma_{WX^3}|$}}
& \textbf{\textit{$|\Gamma_{X^3}| \; |\Gamma_{WX^3}|$}} \\
\hline
        (5,3)  &   [66572 $\;$46415] &    0.39       &   0.39116     &   0.13333     &  0.15255 \\
%\hline
%       (5,3)  &   2000 primes &  0.381        &  0.390195       &  0.1315    &  0.152252 \\
%\hline
        (4,11) &   [-153 $\;$-676]   &    0.37       &   0.37782     &   0.12594     &  0.13979\\
%\hline
        (3,7)  &   [-14 $\;$-127]    &   0.3975      &   0.39783     &   0.14095     &  0.15814 \\
%\hline
        (3,-2) &   [1244$\;$4431]   &  0.3958       &   0.37917     &   0.1275      &  0.15009  \\
%\hline
        (7,11) &   [684 $\;$1493]   &  0.385        &   0.38365     &   0.13        &   0.14771\\
%\hline
        (2,-5) &   [1251 $\;$4364]   &  0.3883       &   0.37615     &  0.12344      &   0.14607\\
\hline
% \multicolumn{4}{l}{ $^{\mathrm{a}}$ Sample of a Table footnote.}
\end{tabular}
\label{tab1}
\end{center}
\end{table}
%=======================================================
%the place to put  Table 4 and Table 5

In both cases we see deterioration of independence.
There is a striking phenomenon for the sets $\Gamma_{X^2}$ and
$\Gamma_{WX^2}$: while the first is quite large, consistent with
Proposition \ref{Proposition24}, the density of the second set is smaller
then the values for generic sequences of Table 1.
Let us note that it has to be smaller then $1/2$ by Theorem
\ref{Theorem31}, yet it does not explain why it is much  smaller
than typical.

In Table 5 we see that the densities are larger than in Table 3,
which is easily explained by the fact that $\Gamma_X\subset \Gamma_{X^3}$.
They cannot exceed $1/2$ because of the constraint of Theorem
\ref{Theorem31}, which applies equally well to $\Gamma_{X^3}$ and $ \Gamma_{WX^3}$.
There is also not much difference between the two. That can be explained
by the fact that $WX^3= (WX)^3$ in $\mathcal G(t)$, since
$W^3 =D^{-1}W$.

%The vague requirement in the Conjecture that $t$ is ``generic'' comes from
%the example of the circular case where the intersections of
%$\Gamma_{X}$ and $ \Gamma_{WX}$ are unusually small, and the conjecture fails.
% Non-primitivity of $t$ most probably also affects the independence.
%In the cubic case the situation is not clear.
%Clearly the sets of prime divisors $\Gamma_{I}$ and $\Gamma_{W}$,
%or $\Gamma_{C}$ and $\Gamma_{V}$, are not independent.
%Also for the torsion elements in the cyclotomic cases
%the sets of divisors do not satisfy the independence conjecture.
%=================================================

%=========================================================================================
%=========================================================================================

\end{document}